# SCHRÖDINGER EQUATIONS WITH SINGULAR POTENTIALS: LINEAR AND NONLINEAR BOUNDARY VALUE PROBLEMS.

MOSHE MARCUS AND PHUOC-TAI NGUYEN

ABSTRACT. Let $\Omega \subset \mathbb{R}^N$ ($N \geq 3$) be a $C^2$ bounded domain and $F \subset \partial\Omega$ be a $C^2$ submanifold of dimension $0 \leq k \leq N-2$. Put $\delta_F(x) = \text{dist}\,(x, F)$, $V = \delta_F^{-2}$ in $\Omega$ and $L_{\gamma V} = \Delta + \gamma V$. Denote by $C_H(V)$ the Hardy constant relative to $V$ in $\Omega$. We study positive solutions of equations (LE) $-L_{\gamma V}u = 0$ and (NE) $-L_{\gamma V}u + f(u) = 0$ in $\Omega$ when $\gamma < C_H(V)$ and $f \in C(\mathbb{R})$ is an odd, monotone increasing function. We establish the existence of a normalized boundary trace for positive solutions of (LE) - first studied in [17] for the case $F = \partial\Omega$ - and employ it to investigate the behavior of subsolutions and super solutions of (LE) at the boundary. Using these results we study boundary value problems for (NE) and derive a-priori estimates. Finally we discuss subcriticality of (NE) at boundary points of $\Omega$ and establish existence and stability results when the data is concentrated on the set of subcritical points.

## CONTENTS

## 1. Introduction

Let $\Omega \subset \mathbb{R}^N$ ($N \geq 3$) be a $C^2$ bounded domain and $F \subset \partial\Omega$ be a $k$ dimensional $C^2$ submanifold, $0 \leq k \leq N-2$. Denote

$$\delta(x) = \text{dist}\,(x, \partial\Omega), \quad \delta_F(x) := \text{dist}\,(x, F)$$

and put

$$V(x) = V_F(x) := \delta_F(x)^{-2}, \quad x \in \Omega,$$
$$L_{\gamma V} = L_{\gamma V}^\Omega := \Delta + \gamma V, \quad \gamma \in \mathbb{R}.$$

In this paper, we study positive solutions of semilinear Schrödinger equation of the form

$$(1.1) \qquad\qquad -L_{\gamma V} u + f(u) = 0 \quad \text{in } \Omega$$

where $f \in C(\mathbb{R})$ is an odd, monotone increasing function. A function $u$ is a solution of equation (1.1) if $u \in L^1_{loc}(\Omega)$, $f(u) \in L^1_{loc}(\Omega)$ and the equation holds in the distribution sense.

To study (1.1) we must first investigate the linear Schrödinger equation

$$(1.2) \qquad\qquad -L_{\gamma V} u = 0 \quad \text{in } \Omega.$$

We say that a function $u \in L^1_{loc}(\Omega)$ is $L_{\gamma V}$ harmonic in $\Omega$ if (1.2) holds in the sense of distributions. If (1.2) holds with "=" replaced by "≤" (resp. "≥") then $u$ is called $L_{\gamma V}$ subharmonic (resp. superharmonic) in $\Omega$.

When $F = \partial\Omega$ and $f(u) = u^q$, equation (1.1) becomes

$$(1.3) \qquad\qquad -\Delta u - \frac{\gamma}{\delta^2} u + u^q = 0 \quad \text{in } \Omega.$$

This equation and its linear counterpart have been studied intensively. Numerous papers deal with the existence and boundary behavior of positive eigenfunctions for the equation

$$(\text{LH}) \qquad\qquad -\Delta u - \frac{\gamma}{\delta^2} u = 0$$

and perturbations of the above. These are closely related to the study of Hardy inequalities. See for instance [6, 7, 8, 9, 10, 15, 18] and the references therein. A sharp estimate of the Green kernel of (LH) in smooth domains was obtained in [12].

Large solutions of (1.3) – i.e., solutions that blow up everywhere on the boundary – have first been investigated in [5]. Assuming that $\gamma < c_H^\Omega \leq 1/4$ ($c_H^\Omega$ = classical Hardy constant) the authors obtained an extension of the Keller–Osserman estimate and investigated the behavior of $L_{\gamma V}$ subharmonic functions at the boundary.

Boundary value problems for (LH) and (1.3), for $\gamma < c_H^\Omega$, have first been discussed in [17]. Here the authors introduced a notion of *normalized boundary trace*, given by a bounded Borel measure on $\partial\Omega$, by which the positive solutions of (LH) could be fully classified. This classification remains valid for (1.3) as long as we restrict ourselves to positive *moderate* solutions, i.e., solutions dominated by an $L_{\gamma V}$ harmonic function.

In [13] boundary value problems have been studied in the framework of weighted $H^1$ spaces. Boundary value problems for the nonlinear equation (1.3) were further investigated in [16] where it was shown that the results of [17] remain valid for all $\gamma < 1/4$ and all bounded $C^2$ domains $\Omega$, even if $c_H^\Omega < 1/4$.

In the present paper we study boundary problems for equations (1.2) and (1.1). As in [17] our approach is based on the introduction of an appropriate normalized boundary trace. The determination of this trace requires precise estimates of the Green and Martin kernels (recently derived in [14]) as well as precise estimates of the ground states. The present study is considerably more delicate than in [17] because of the nonuniform behavior of the potential at the boundary.

Let $C_H(V)$ be *the Hardy constant relative to $V$* in $\Omega$, namely

$$(1.4) \qquad C_H(V) = \inf_{\phi \in C_c^1(\Omega)} \frac{\int_\Omega |\nabla \phi|^2 \, dx}{\int_\Omega \phi^2 V \, dx}.$$

Since $V = V_F \leq \delta^{-2}$, it follows that $C_H(V) \geq c_H^\Omega > 0$. On the other hand, by [11, Corollary 1.3], $C_H(V) \leq \frac{(N-k)^2}{4}$.

It is well known that, if $\gamma < C_H(V)$, there exists a positive first eigenvalue $\lambda_{\gamma V}$ of $-L_{\gamma V}$ in $\Omega$ with corresponding positive eigenfunction $\varphi_{\gamma V} \in H_0^1(\Omega)$. The eigenfunction is normalized by $\varphi_{\gamma V}(x_0) = 1$ where $x_0$ is a fixed reference point in $\Omega$. In particular $\varphi_0$ denotes the first eigenfunction of the Laplacian normalized by $\varphi_0(x_0) = 1$.

The following result is due to Pinchover [20].

*Suppose that $V \in C(\Omega)$ satisfies*

$$(A) \qquad \text{(i) } 0 < V \leq c\delta(x)^{-2}, \quad \text{(ii) } \gamma < C_H(V).$$

*Then there exists $\epsilon > 0$ such that the operator $-(L_{\gamma V} + \epsilon \delta(x)^{-2})$ has a positive supersolution.*

A proof of this result is provided in [14].

By Ancona [2], it follows that, assuming (A), $L_{\gamma V}$ is weakly coercive and consequently the results of [2] and [3] apply to $V$. These include the representation theorem for positive $L_{\gamma V}$ harmonic functions and the Boundary Harnack Principle (briefly BHP).

Denote by $G_{\gamma V}^\Omega$ (respectively $K_{\gamma V}^\Omega$) the Green (respectively Martin) kernel of $-L_{\gamma V}$ in $\Omega$. The superscript will generally be dropped when the underlying domain is $\Omega$

It was pointed out to us by Alano Ancona that – assuming (A) – given $x_0 \in \Omega$ and a neighborhood $Q$ of $\partial\Omega$ such that $x_0 \notin \bar{Q}$ there exists a constant $C = C(x_0, Q)$ such that

$$(1.5) \qquad \frac{1}{C} \varphi_{\gamma V}(x) \le G_{\gamma V}(x, x_0) \le C \varphi_{\gamma V}(x) \quad \forall x \in Q \cap \Omega.$$

See [14] for more details.

Using these facts, the first author [14] obtained sharp, two sided estimates of the Green and Martin kernels for operators $L_{\gamma V}$ satisfying (A). These, together with sharp two sided estimates of $\varphi_{\gamma V}$ based on the work of Fall and Mahmoudi [11] (see Section 2 below) lie at the core of the present paper.

Denote by $\mathfrak{M}(\partial \Omega)$ the space of finite Borel measures on $\partial \Omega$ and by $\mathfrak{M}^+(\partial \Omega)$ the positive cone of $\mathfrak{M}(\partial \Omega)$. Denote by $\mathfrak{M}(\Omega, \phi) - \phi$ a positive Borel function in $\Omega$ – the space of Radon measures $\tau$ on $\Omega$ satisfying $\int_{\Omega} \phi \, d|\tau| < \infty$ and by $\mathfrak{M}^+(\Omega, \phi)$ the positive cone of this space.

Put

$$(1.6) \qquad \begin{aligned} \mathbb{K}_{\gamma V}[\nu] &:= \int_{\partial \Omega} K_{\gamma V}(\cdot, y) d\nu(y), \quad \nu \in \mathfrak{M}(\partial \Omega), \\ \mathbb{G}_{\gamma V}[\tau] &:= \int_{\Omega} G_{\gamma V}(\cdot, y) d\tau(y), \quad \tau \in \mathfrak{M}^+(\Omega). \end{aligned}$$

Inequality (1.5) implies:

$\mathbb{G}_{\gamma V}[\tau]$ is finite if and only if $\tau \in \mathfrak{M}^+(\Omega, \varphi_{\gamma V})$, i.e.,

$$\int_{\Omega} \varphi_{\gamma V} d\tau < \infty.$$

For $\beta > 0$, denote

$$(1.7) \qquad \begin{aligned} \Omega_\beta &:= \{x \in \Omega : \delta(x) < \beta\}, \ D_\beta := \{x \in \Omega : \delta(x) > \beta\}, \\ \Sigma_\beta &:= \{x \in \Omega : \delta(x) = \beta\}. \end{aligned}$$

In what follows the notation $f \sim g$ in a domain $D$ means that there exists two positive constants $c_1, c_2$ such that $c_1 f \le g \le c_2 f$ in $D$. The constants $c_1, c_2$ are called *similarity constants*.

Our first result is the key to the determination of the correct normalization for the definition of boundary trace.

Due to two-sided estimates on $K_{\gamma V}$ [14], the following key estimates holds when $V = V_F$, (see Theorem 3.1 and Corollary **??**)

**Theorem 1.1.** *Assume*

$$(1.8) \qquad \gamma < \min\{C_H(V), \frac{2(N-k)-1}{4}\}.$$

*Then, for $V = V_F$,*

$$(1.9)$$
$$\int_{\Sigma_\beta} \frac{\varphi_{\gamma V}(x)}{\varphi_0(x)} \mathbb{K}_{\gamma V}[\nu](x) dS_x \sim \|\nu\|_{\mathfrak{M}(\partial \Omega)}, \quad \forall \nu \in \mathfrak{M}^+(\partial \Omega), \quad \forall \beta \ small.$$

*The similarity constants depend only on $N, \Omega, F, \gamma$.*

This leads to the following

**Definition 1.2.** Assume (1.8). We say that a function $u$ possesses *normalized boundary trace* $\nu \in \mathfrak{M}(\partial\Omega)$ if

$$(1.10) \qquad \lim_{\beta \to 0} \int_{\Sigma_\beta} \frac{\varphi_{\gamma V}(x)}{\varphi_0(x)} |u(x) - \mathbb{K}_{\gamma V}[\nu](x)| dS_x = 0.$$

The normalized boundary trace of $u$ is denoted by $\mathrm{tr}^*(u)$.

**Remark 1.3.** (i) In the case of $F = \partial\Omega$ this definition reduces to the definition of normalized boundary trace introduced in [17, Definition 1.2].

(ii) The notion of normalized boundary trace is well defined. Indeed, suppose that there exist two measures $\nu$ and $\nu'$ in $\mathfrak{M}(\partial\Omega)$ satisfying (1.10). Put $v := (\mathbb{K}_{\gamma V}[\nu - \nu'])_+$ then $v$ is a nonnegative $L_{\gamma V}$ subharmonic function with $\mathrm{tr}^*(v) = 0$. Then by Proposition 4.4, $v = 0$ and hence $\mathbb{K}_{\gamma V}[\nu - \nu'] \leq 0$ in $\Omega$. By changing the role of $\nu$ and $\nu'$, we deduce that $\nu = \nu'$.

Next we consider the linear boundary value problem

$$(1.11) \qquad \begin{cases} -L_{\gamma V} u = \tau & \text{in } \Omega, \\ \mathrm{tr}^*(u) = \nu, \end{cases}$$

where $\tau \in \mathfrak{M}(\Omega, \varphi_{\gamma V})$ and $\nu \in \mathfrak{M}(\partial\Omega)$. A function $u \in L^1_{loc}(\Omega)$ is a solution of (1.11) if $u$ satisfies the equation in (1.11) in the sense of distribution and $\mathrm{tr}^*(u) = \nu$.

Following are our main results on problem (1.11) with $V = V_F$.

**Theorem 1.4.** I. *If* $\tau \in \mathfrak{M}(\Omega, \varphi_{\gamma V})$ *then* $\mathrm{tr}^*(\mathbb{G}_{\gamma V}[\tau]) = 0$. *Thus* $\mathbb{G}_{\gamma V}[\tau]$ *is a solution of* (1.11) *with* $\nu = 0$.

II. *Let* $u$ *be a positive* $L_{\gamma V}$ *subharmonic function. If* $u$ *is dominated by an* $L_{\gamma V}$ *superharmonic function then* $\lambda := L_{\gamma V} u \in \mathfrak{M}^+(\Omega, \varphi_{\gamma V})$ *and* $u$ *has a normalized boundary trace, say* $\nu$. *In this case*

$$(1.12) \qquad u + \mathbb{G}_{\gamma V}[\lambda] = \mathbb{K}_{\gamma V}[\nu].$$

*In particular*

$$(1.13) \qquad u \leq \mathbb{K}_{\gamma V}[\nu], \quad \nu = \mathrm{tr}^*(u).$$

III. *Let* $u$ *be a positive* $L_{\gamma V}$ *superharmonic function. Then there exist* $\nu \in \mathfrak{M}^+(\partial\Omega)$ *and* $\tau \in \mathfrak{M}^+(\Omega, \varphi_{\gamma V})$ *such that*

$$(1.14) \qquad u = \mathbb{G}_{\gamma V}[\tau] + \mathbb{K}_{\gamma V}[\nu].$$

*In particular,*

$$(1.15) \qquad u \geq \mathbb{K}_{\gamma V}[\nu], \quad \nu = \mathrm{tr}^*(u)$$

*and $u$ is an $L_{\gamma V}$ potential (i.e., $u$ does not dominate any positive $L_{\gamma V}$ harmonic function) if and only if $u$ is the Green potential of a measure in $\mathfrak{M}^+(\Omega, \varphi_{\gamma V})$.*

IV. *For every $\nu \in \mathfrak{M}(\partial\Omega)$ and $\tau \in \mathfrak{M}(\Omega, \varphi_{\gamma V})$, problem (1.11) has a unique solution. The solution is given by (1.14).*

Theorem 1.4 plays an important role in the study of nonlinear boundary value problems of the form

$$(1.16) \qquad \begin{cases} -L_{\gamma V} u + f(u) = 0 & \text{in } \Omega, \\ \quad\quad\quad \operatorname{tr}^*(u) = \nu. \end{cases}$$

We shall consider functions $f$ as in (1.1).

**Definition 1.5.** (i) A positive solution of (1.1) is *$L_{\gamma V}$ moderate* if it is dominated by an $L_{\gamma V}$ harmonic function.

(ii) A function $u$ is a *(weak) solution of* (1.16) if $u$ is a solution of (1.1) and $\operatorname{tr}^*(u) = \nu$.

The next results describe the main properties of moderate solutions of (1.1).

**Theorem 1.6.** *Let $u$ be a positive solution of* (1.1). *Then the following statements are equivalent:*

*(i) $u$ is $L_{\gamma V}$ moderate.*

*(ii) $u$ admits a normalized boundary trace $\nu \in \mathfrak{M}^+(\partial\Omega)$. In other words, $u$ is a solution of* (1.16).

*(iii) $f(u) \in L^1(\Omega, \varphi_{\gamma V})$.*

*Finally, if $u$ is a positive solution of* (1.16) *then*

$$(1.17) \qquad u + \mathbb{G}_{\gamma V}[f(u)] = \mathbb{K}_{\gamma V}[\nu].$$

Existence and uniqueness results as well as a priori estimates are presented in the next theorems.

**Theorem 1.7.** *(i) For every $\nu \in \mathfrak{M}^+(\partial\Omega)$ there exists at most one positive solution of* (1.16).

*(ii) Let $\nu_i \in \mathfrak{M}^+(\partial\Omega)$, $i = 1, 2$ and suppose that there exist corresponding solutions $u_i$ of* (1.16) *with $\nu$ replaced by $\nu_i$. If $\nu_1 \leq \nu_2$ then $u_1 \leq u_2$ in $\Omega$.*

*(iii) If $u$ is a positive solution of* (1.16) *then*

$$(1.18) \qquad \|u + \mathbb{G}_{\gamma V}[f(u)]\|_{L^1(\Omega, \frac{\varphi_{\gamma V}}{\varphi_0})} \sim \|\nu\|_{\mathfrak{M}(\partial\Omega)}$$

*and*

$$(1.19) \qquad \|f(u)\|_{L^1(\Omega, \varphi_{\gamma V})} + \|u\|_{L^1(\Omega, \frac{\varphi_{\gamma V}}{\varphi_0})} \sim \|\nu\|_{\mathfrak{M}(\partial\Omega)}.$$

*The similarity constants depend only on $\Omega$, $\gamma V$ and $f$.*

**Definition 1.8.** A measure $\nu \in \mathfrak{M}^+(\partial\Omega)$ is called a *good measure* relative to $(f, L_{\gamma V})$ if (1.16) has a solution. The family of good measures is denoted by $\mathcal{M}_{\gamma V}(f)$.

**Proposition 1.9.** *For every $\nu \in \mathfrak{M}^+(\partial\Omega)$,*

$$(1.20) \qquad f(\mathbb{K}_{\gamma V}[\nu]) \in L^1(\Omega, \varphi_{\gamma V}) \Longrightarrow \nu \in \mathcal{M}_{\gamma V}(f).$$

**Remark 1.10.** Condition (1.20) is sufficient but *not necessary* for a measure to be good. For instance if $\nu$ is the limit of an increasing sequence of good measures then $\nu \in \mathcal{M}_{\gamma V}(f)$. In particular, every non-negative function $h \in L^1(\partial\Omega)$ is a good measure because it is the limit of an increasing sequence of bounded functions. However, if $f(t) \geq t^p$ for some $p > 1$ and all $t > 1$, it is easy to construct functions in $L^1(\partial\Omega)$ which do not satisfy (1.20).

Next we define the notion of subcriticality of (1.1) at a point $y \in \partial\Omega$ and investigate in more detail problem (1.16) when $\nu$ is supported in the set of subcritical points.

**Definition 1.11.** A point $y \in \partial\Omega$ is $(f, L_{\gamma V})$ *subcritical* if the problem

$$(1.21) \qquad -L_{\gamma V}u + f(u) = 0 \quad \Omega, \quad \mathrm{tr}^*(u) = k\delta_y$$

has a solution for every $k > 0$. Here $\delta_y$ denotes the Dirac measure concentrated at $y$. We denote by $SC_{\gamma V}(f)$ the set of $(f, L_{\gamma V})$ subcritical points.

**Definition 1.12.** Denote by $\mathcal{F}_0$ the family of functions $f \in C(\mathbb{R})$ such that $f$ is odd, monotone increasing and $f(t)/t$ is non-decreasing and tends to infinity as $t \to \infty$. Denote by $\mathcal{F}(a^*)$ the family of functions $f \in \mathcal{F}_0$ that satisfy the following condition:

There exists a constant $a^* > 0$ such that, for every smooth domain $D \subset \Omega$, if $u$ is a positive function satisfying

$$(1.22) \qquad\qquad -\Delta u + f(u) \leq 0 \quad \text{in} \quad \Omega$$

then

$$(1.23) \qquad f(u(x)) \leq a^* \mathrm{dist}\,(x, \partial D)^{-2} u(x) \quad \forall x \in D.$$

**Theorem 1.13.** *Assume $f \in \mathcal{F}(a^*)$. If $\nu \in \mathfrak{M}^+(\partial\Omega)$ is concentrated on $SC_{\gamma V}(f)$ then $\nu \in \mathcal{M}_{\gamma V}(f)$.*

**Theorem 1.14.** *Assume $f \in \mathcal{F}(a^*)$. Let $\nu \in \mathcal{M}_{\gamma V}(f)$ and let $u$ be the corresponding solution of* (1.16). *Then*

$$(1.24) \qquad \lim_{x \to y} \frac{u(x)}{\mathbb{K}_{\gamma V}[\nu](x)} = 1 \quad \textit{non tangentially, $\nu$-a.e. on $\partial\Omega$.}$$

**Definition 1.15.** We say that $f$ satisfies the $\Delta_2$ condition if there exists a positive constant $c$ such that

$$f(a + b) \leq c(f(a) + f(b)) \quad \forall a > 0, b > 0.$$

A necessary and sufficient condition for a boundary point to be $(f, L_{\gamma V})$ subcritical is given in the next resut.

**Theorem 1.16.** *Let $f \in \mathcal{F}(a^*)$ such that $f$ satisfies the $\Delta_2$ condition. A point $y$ is $(f, L_{\gamma V})$ subcritical, i.e. $y \in SC_{\gamma V}(f)$, if and only if $f(K_{\gamma V}(\cdot, y)) \in L^1(\Omega, \varphi_{\gamma V})$.*

For $\gamma < C_H(V)$, let $\alpha$ be the largest root of equation $\alpha^2 + (N-k)\alpha + \gamma = 0$, namely

$$(1.25) \qquad \alpha = \alpha_{k,\gamma} := \frac{k - N + \sqrt{(N-k)^2 - 4\gamma}}{2}.$$

Put

$$(1.26) \qquad q_F^* := \frac{N+1+\alpha}{N-1+\alpha} \quad \text{and} \quad q^* := \frac{N+1}{N-1}.$$

Existence and stability results are provided in the following theorems.

**Theorem 1.17.** *Let $f(t) = |t|^q \mathrm{sign}\, t$, $q > 1$.*

*(a) Assume $\alpha > -1$. A point $y \in F$ is $(f, L_{\gamma V})$ subcritical if and only if $q \in (1, q_F^*)$.*

*(b) A point $y \in \partial\Omega \setminus F$ is $(f, L_{\gamma V})$ subcritical if and only if $q \in (1, q^*)$.*

**Theorem 1.18.** *Let $f \in C(\mathbb{R})$ be an odd, monotone increasing function. Assume that for some $p \in (1, \infty)$,*

$$(1.27) \qquad \int_0^1 f(ks^{-1/p})ds < \infty \quad \forall k > 0.$$

*(i) If $p = q_F^*$ then $F \subset SC_{\gamma V}(f)$ and every measure $\nu \in \mathfrak{M}^+(\partial\Omega)$ such that $\mathrm{supp}\,\nu \subset F$ is a good measure.*

*If $\{\nu_n\}$ is a bounded sequence of positive measures supported in $F$ then $\{f(\mathbb{K}_{\gamma V}[\nu_n])\}$ is uniformly absolutely continuous in $L^1(\Omega, \varphi_{\gamma V})$.*

*(ii) If $p = q^*$ then $\partial\Omega \setminus F \subset SC_{\gamma V}(f)$ and every measure $\nu \in \mathfrak{M}^+(\partial\Omega)$ such that $\mathrm{supp}\,\nu \subset \partial\Omega \setminus F$ is a good measure.*

*If $\{\nu_n\}$ is a bounded sequence of positive measures supported in a fixed compact subset of $\partial\Omega \setminus F$ then $\{f(\mathbb{K}_{\gamma V}[\nu_n])\}$ is uniformly absolutely continuous in $L^1(\Omega, \varphi_{\gamma V})$.*

**Theorem 1.19.** *Let $f \in \mathcal{F}(a^*)$. Assume that $f$ satisfies the $\Delta_2$ condition and, for some $p \in (1, \infty)$,*

$$(1.28) \qquad \int_0^1 f(s^{-1/p})ds < \infty.$$

*Let $\{\nu_n\} \subset \mathfrak{M}^+(\partial\Omega)$ and assume that $\nu_n \rightharpoonup \nu$, i.e. the sequence converges to $\nu$ weakly relative to $C(\partial\Omega)$.*

*(i) Assume that the measures $\nu_n$ are supported in $F$ and that $p = q_F^*$. Let $u_n$ be the solution of problem (1.16) with $\nu = \nu_n$ and let $u$ be the solution of (1.16) where $\nu$ is the weak limit of the sequence. Then*

$$(1.29) \qquad \begin{aligned} u_n &\to u \quad \text{locally uniformly in } \Omega, \\ f(u_n) &\to f(u) \quad \text{in } L^1(\Omega, \varphi_{\gamma V}). \end{aligned}$$

*(ii) Assume that the measures $\nu_n$ are supported in a fixed compact subset of $\partial\Omega \setminus F$ and that $p = q^*$. Let $u_n$ be the solution of problem (1.16) with $\nu = \nu_n$ and let $u$ be the solution of (1.16) when $\nu$ is the weak limit of the sequence. Then (1.29) holds.*

The paper is organized as follows. In Section 2, we establish estimates of the first eigenfunction of $-L_{\gamma V}$ by adapting an argument of [11]. This result together with the estimates of the Green and Martin kernels of [14] are crucial in the study of the notion of normalized boundary trace in Section 3. The linear boundary value problem is discussed in Section 4. In Section 5, we investigate moderate solutions of nonlinear boundary value problem. Finally, we obtain existence and stability results in Sections 6 and 7.

## 2. Estimates on eigenfunction of $-L_{\gamma V}$

Let $\Omega \subset \mathbb{R}^N$ ($N \geq 3$) be a $C^2$ bounded domain and $F \subset \partial\Omega$ be a $C^2$ submanifold with dimension $0 \leq k \leq N - 2$. Recall that $\delta(x) = \text{dist}(x, \partial\Omega)$ and $\delta_F(x) = \text{dist}(x, F)$. In this section, we establish two-sided estimates of the first eigenfunction of $-L_{\gamma V}$ in $\Omega$ when $V = V_F$.

In the sequel the notation $f \lesssim g$ (resp. $f \gtrsim g$) in a domain $D$ means that there exists a positive constant $c$ such that $f \leq c\,g$ (resp. $f \geq c\,g$) in $D$.

For $x \in \Omega$, denote by $\sigma(x)$ the projection of $x$ on $\partial\Omega$. Let $\text{dist}^{\partial\Omega}(\cdot, F)$ be the geodesic distance to $F$ on $\partial\Omega$. Put

$$(2.1) \qquad \tilde{\delta}_F(x) := \sqrt{\delta(x)^2 + (\text{dist}^{\partial\Omega}(\sigma(x), F))^2}.$$

Note that the distances $\delta_F$ and $\tilde{\delta}_F$ are equivalent (see [11, Lemma 2.1]). For $M, a \in \mathbb{R}$, define

$$(2.2) \qquad W_{a,M}(x) := X_a(\tilde{\delta}_F(x)) e^{M\delta(x)} \delta(x) \tilde{\delta}_F(x)^\alpha$$

where $\alpha$ is given in (1.25) and

$$(2.3) \qquad X_a(t) := (-\ln(t))^a, \quad t \in (0, 1/2).$$

For $\epsilon > 0$, put

$$\Omega_\epsilon(F) := \{x \in \Omega : \delta_F(x) < \epsilon\}.$$

Since $F$ is a $C^2$-submanifold of $\partial\Omega$ there exists a positive constant $\epsilon_1$ depending on $F, \Omega$ such that $\delta_F \in C^2(\Omega_{\epsilon_1}(F))$.

**Lemma 2.1.** *Let $\lambda \in \mathbb{R}$. There exist constants $M_1 = M_1(F, \Omega) > 0$, $M_2 = M_2(F, \Omega) < 0$ and $\epsilon_2 = \epsilon_2(\lambda, \gamma, N, F, \Omega) > 0$ such that the functions*

$$U_1 := W_{-1,M_1} + W_{0,M_1} \quad and \quad U_2 := W_{0,M_2} - W_{-1,M_1} > 0 \ in \ \Omega_{\epsilon_2}(F)$$

*satisfy*

$$(2.4) \qquad (-L_{\gamma V} + \lambda)U_1 \le 0 \le (-L_{\gamma V} + \lambda)U_2 \quad in \ \Omega_{\epsilon_2}(F).$$

*Proof.* One can proceed as in the proof of [11, Lemma 2.3, Lemma 2.4]. The main difference is that $\alpha$ defined in (1.25) does not depend on $\bar{\delta}_F$, which simplifies the proof in the present case. $\qquad\square$

**Lemma 2.2.** *Let $\gamma < C_H(V)$. Then there exists a positive constant $\epsilon_3 = \epsilon_3(\gamma, N, F, \Omega) > 0$ such that*

$$(2.5) \qquad \varphi_{\gamma V}(x) \sim \delta(x)\delta_F(x)^\alpha, \quad \forall x \in \Omega_{\epsilon_3}(F),$$

*where $\alpha$ is given in (1.25). Here the similarity constants depend on $N, \gamma, F, \Omega$.*

*Proof.* By local improved Hardy inequality [11, Lemma 3.1], there exist $\epsilon_3 = \epsilon_3(\gamma, N, F, \Omega)$ and $c = c(\gamma, N, F, \Omega) > 0$ such that

$$(2.6)$$
$$\int_{\Omega_{\epsilon_3}(F)} \left(|\nabla\phi|^2 - \gamma V\phi^2\right) dx \ge c \int_{\Omega_{\epsilon_3}(F)} \frac{\phi^2}{\delta_F^2 |\ln(\delta_F)|^2} dx, \ \forall \phi \in H_0^1(\Omega_{\epsilon_3}(F)).$$

Moreover, we can assume that $\epsilon_3$ small enough such that

$$\delta_F(x)^2 |\ln(\delta_F(x))|^2 \le c(2\lambda_{\gamma V})^{-1}, \quad \forall x \in \Omega_{\epsilon_3}(F).$$

Consequently,

$$(2.7)$$
$$\int_{\Omega_{\epsilon_3}(F)} \left(|\nabla\phi|^2 - \gamma V\phi^2\right) dx \ge 2\lambda_{\gamma V} \int_{\Omega_{\epsilon_3}(F)} \phi^2 dx, \quad \forall \phi \in H_0^1(\Omega_{\epsilon_3}(F)).$$

The eigenfunction $\varphi_{\gamma V}$ satisfies

$$-\Delta\varphi_{\gamma V} - (\gamma V + \lambda_{\gamma V})\varphi_{\gamma V} = 0 \quad in \ \Omega.$$

Since

$$|\gamma V + \lambda_{\gamma V}| \le c(\gamma, N, F, \Omega) \quad in \ \overline{\Omega \setminus \Omega_{\frac{\epsilon_3}{2}}(F)},$$

it follows from standard regularity theory that

$$\varphi_{\gamma V} \le c(\gamma, N, F, \Omega) \quad in \ \overline{\Omega \setminus \Omega_{\frac{\epsilon_3}{2}}(F)}.$$

Therefore, one can choose another constant $c = c(\gamma, N, F, \Omega) > 0$ such that $\varphi_{\gamma V} \le c\, U_2$ in $\overline{\Omega \cap \partial\Omega_{\epsilon_3}(F)}$ where $U_2$ is given in Lemma 2.1. Put $\tilde{u} := \varphi_{\gamma V} - c\, U_2$ then $\tilde{u}_+ \in H_0^1(\Omega_{\epsilon_3}(F))$. By Lemma 2.1,

$$-L_{\gamma V}\tilde{u} \le \lambda_{\gamma V}\tilde{u} \quad in \ \Omega_{\epsilon_3}(F).$$

Multiplying by $\tilde{u}_+$ and integrating by parts yield

$$(2.8) \qquad \int_{\Omega_{\epsilon_3}(F)} \left( |\nabla \tilde{u}_+|^2 - \gamma V \tilde{u}_+^2 \right) dx \leq \lambda_{\gamma V} \int_{\Omega_{\epsilon_3}(F)} \tilde{u}_+^2 dx.$$

By combining (2.8) and (2.7) with $\phi$ replaced by $\tilde{u}_+$, we obtain

$$2\lambda_{\gamma V} \int_{\Omega_{\epsilon_3}(F)} \tilde{u}_+^2 dx \leq \lambda_{\gamma V} \int_{\Omega_{\epsilon_3}(F)} \tilde{u}_+^2 dx.$$

This implies $\tilde{u}_+ \equiv 0$ in $\Omega_{\epsilon_3}(F)$. Therefore

$$\varphi_{\gamma V}(x) \leq c\, U_2(x) \leq c\, W_{0,M_2}(x) \leq c\, \delta(x) \delta_F(x)^\alpha \quad \forall x \in \Omega_{\epsilon_3}(F).$$

Similarly, there exists another constant $c = c(\gamma, N, F, \Omega) > 0$ such that $c\, U_1 \leq \varphi_{\gamma V}$ in $\Omega \cap \partial\Omega_{\epsilon_3}(F)$, where $U_1$ is given in Lemma 2.1. Put $\tilde{v} := c\, U_1 - \varphi_{\gamma V}$ then $\tilde{v}_+ \in H_0^1(\Omega_{\epsilon_3}(F))$ since $U_1 = 0$ on $\partial\Omega \setminus \partial\Omega_{\epsilon_3}(F)$. Moreover,

$$L_{\gamma V} \tilde{v} \leq \lambda_{\gamma V} \tilde{v} \quad \text{in } \Omega_{\epsilon_3}(F).$$

Again, by using (2.7) with $\phi$ replaced by $\tilde{v}_+$, we deduce that $\tilde{v}_+ = 0$ in $\Omega_{\epsilon_3}(F)$, i.e. $\varphi_{\gamma V} \geq c\, U_1$ in $\Omega_{\epsilon_3}(F)$. This leads to

$$\varphi_{\gamma V}(x) \geq c\, U_1(x) \geq c\, W_{0,M_1}(x) \geq c\, \delta(x) \delta_F(x)^\alpha \quad \forall x \in \Omega_{\epsilon_3}(F).$$

Combining the above inequalities, we obtain (2.5). $\qquad \square$

## 3. Normalized boundary trace

Let $(r_0, \kappa)$ ($\kappa \geq 1$) be the $C^2$ characteristic of $\Omega$ (see [19, Page 1] for more detail). We denote by $G_0$ and $K_0$ respectively the Green kernel and Poisson kernel of $-\Delta$ in $\Omega$. We also denote

$$(3.1) \quad A(x,y) := \{z \in \Omega : \frac{1}{2} \hat{r}(x,y) \leq \delta(z) \leq 2\hat{r}(x,y)\} \cap B_{4|x-y|}(\frac{x+y}{2}).$$

where $\hat{r}(x,y) = \max\{|x-y|, \delta(x), \delta(y)\}$. Finally, throughout this section we assume that

$$(3.2) \qquad \gamma < \min\{C_H(V), \frac{2(N-k)-1}{4}\}.$$

**Theorem 3.1.** *There exists a positive constant $\epsilon_0$ such that*

$$(3.3) \qquad \int_{\Sigma_\beta} \frac{\varphi_{\gamma V}(x)}{\varphi_0(x)} K_{\gamma V}(x,y) dS_x \sim 1 \quad \forall y \in \partial\Omega, \ \forall \beta \in (0, \epsilon_0),$$

*with similarity constants independent of $y$ and $\beta$. Here $\varphi_0$ is the first eigenfunction of $-\Delta$ in $\Omega$ normalized by $\varphi_0(x_0) = 1$.*

*Proof.* Let $y \in \partial\Omega$ and $x \in \Omega$ be such that $\delta(x) < r_0/10\kappa$ and $|x-y| < r_0$. In [14, Theorem 1.5], we choose $x_y \in A(x,y)$ such that

$$(3.4) \qquad |x-y| < \delta(x_y) = |x_y - y| < \frac{3}{2}|x-y|,$$
$$\delta_F(x_y) \sim \delta(x_y) + \delta_F(y) \sim |x-y| + \delta_F(y).$$

These estimates hold, for instance, if $x_y = y + |x - y|\mathbf{n}_y$. In this case

$$\max(\delta_F(y), \delta(x_y)) \lesssim \delta_F(x_y).$$

As $|\delta_F(x) - \delta_F(y)| \leq |x - y|$ it follows that

$$\delta_F(x) \lesssim \delta_F(x_y) \lesssim \delta_F(x) + |x - y|.$$

On the other hand $|x - y| \sim \delta(x_y) \leq \delta_F(x_y)$. Hence

$$(3.5) \qquad \delta_F(x_y) \sim \delta_F(x) + |x - y|.$$

Let $0 < \epsilon_0 < \frac{1}{4}\min(\epsilon_1, \epsilon_2, \epsilon_3, \frac{r_0}{10\kappa})$ ($\epsilon_i$ as in Section 2).

For any $x \in \Omega_{\epsilon_0}(F)$ and $y \in \partial\Omega$ such that $|x - y| < \epsilon_0$, let $x_y$ be a point satisfying (3.4), (3.5). Recall that $\alpha$ is given in (1.25). Then, by [14, Theorem 1.5] and Lemma 2.2 we obtain

$$
\begin{aligned}
(3.6) \qquad \frac{\varphi_{\gamma V}(x)}{\varphi_0(x)} K_{\gamma V}(x, y) &\sim \frac{\varphi_{\gamma V}(x)^2}{\varphi_{\gamma V}(x_y)^2 \varphi_0(x)} |x - y|^{2-N} \\
&\sim \frac{\delta(x)\delta_F^2(x)^{2\alpha}}{\delta(x_y)^2 \delta_F(x_y)^{2\alpha}} |x - y|^{2-N} \\
&\sim \left(\frac{\delta_F(x)}{\delta_F(x_y)}\right)^{2\alpha} \delta(x)|x - y|^{-N} \\
&\sim \left(\frac{\delta_F(x)}{\delta_F(x) + |x - y|}\right)^{2\alpha} \delta(x)|x - y|^{-N} \\
&\sim \left(1 + \frac{|x - y|}{\delta_F(x)}\right)^{-2\alpha} \delta(x)|x - y|^{-N}.
\end{aligned}
$$

For any $x \in \Omega_{\epsilon_0}(F)$ and $y \in \partial\Omega$ such that $|x - y| \geq \epsilon_0$,

$$(3.7) \qquad \frac{\varphi_{\gamma V}(x)}{\varphi_0(x)} K_{\gamma V}(x, y) \sim \frac{\varphi_{\gamma V}(x)^2}{\varphi_0(x)} \sim \delta(x)\delta_F(x)^{2\alpha}.$$

For $y \in \partial\Omega$ and $\beta \in (0, \epsilon_0)$ put

$$I(y) := \int_{\Sigma_\beta} \frac{\varphi_{\gamma V}(x)}{\varphi_0(x)} K_{\gamma V}(x, y) dS_x,$$

$$I_1(y) := \int_{\Sigma_\beta \cap \Omega_{\epsilon_0}(F)} \frac{\varphi_{\gamma V}(x)}{\varphi_0(x)} K_{\gamma V}(x, y) dS_x,$$

$$I_2(y) := \int_{\Sigma_\beta \setminus \Omega_{\epsilon_0}(F)} \frac{\varphi_{\gamma V}(x)}{\varphi_0(x)} K_{\gamma V}(x, y) dS_x.$$

*Case 1:* $-\frac{1}{2} < \alpha \leq 0$. Note that

$$-\frac{1}{2} < \alpha \iff \gamma < \frac{2(N-k)-1}{4}.$$

Put $\Sigma_\beta \cap \Omega_{\epsilon_0}(F) =: \Sigma_{\beta,\epsilon_0}(F)$. By (3.6) and since $\alpha \leq 0$, if $\beta < \epsilon_0/2$,

$$
\begin{aligned}
I_1(y) &\gtrsim \int_{\Sigma_{\beta,\epsilon_0}(F) \cap B_{\epsilon_0}(y)} \delta(x) |x-y|^{-N} dS_x \\
&\sim \int_{\Sigma_{\beta,\epsilon_0}(F) \cap B_{\epsilon_0}(y)} K_0(x,y) dS_x.
\end{aligned}
\tag{3.8}
$$

By (3.6), (3.7), the inequality $\delta(x) \leq \delta_F(x)$ and since $\alpha \leq 0$, we obtain

$$
\begin{aligned}
I_1(y) &\lesssim \int_{\Sigma_{\beta,\epsilon_0}(F) \cap B_{\epsilon_0}(y)} \left( \frac{\delta(x)}{\delta(x)+|x-y|} \right)^{2\alpha} \delta(x)|x-y|^{-N} dS_x \\
&\quad + \int_{\Sigma_{\beta,\epsilon_0}(F) \setminus B_{\epsilon_0}(y)} \delta(x)\delta_F(x)^{2\alpha} dS_x \\
&\lesssim \beta^{2\alpha+1} \int_{\Sigma_{\beta,\epsilon_0}(F) \cap B_{\epsilon_0}(y)} (\beta+|x-y|)^{-2\alpha}|x-y|^{-N} dS_x \\
&\quad + \beta^{2\alpha+1} \int_{\Sigma_{\beta,\epsilon_0}(F) \setminus B_{\epsilon_0}(y)} dS_x =: I_{1,1}(y) + I_{1,2}(y).
\end{aligned}
\tag{3.9}
$$

Clearly, if $\alpha > -\frac{1}{2}$ then $I_{1,2}(y)$ is bounded by a constant depending on $\epsilon_0$.

We estimate $I_{1,1}(y)$. Let $y'$ be the unique point on $\Sigma_\beta$ such that $|y-y'| = \beta$. Let $\epsilon_0 > 0$ be sufficiently small (depending only on the geometry of $\partial\Omega$ and $F$) so that,

$$
\frac{1}{2}(|x-y'|+\beta) \leq |x-y| \leq 2(|x-y'|+\beta) \quad \forall x \in \Sigma_\beta \cap B_{\epsilon_0}(y).
$$

Using this inequality we obtain,

$$
\begin{aligned}
I_{1,1}(y) &\lesssim \beta^{2\alpha+1} \int_{|x-y'|<\epsilon_0} (|x-y'|+\beta)^{-2\alpha-N} dS_x \\
&\lesssim \beta^{2\alpha+1} \int_0^{\epsilon_0} (t+\beta)^{-2\alpha-2} dt \\
&\lesssim \beta^{2\alpha+1} \int_\beta^{\beta+\epsilon_0} t^{-2\alpha-2} dt.
\end{aligned}
\tag{3.10}
$$

If $\alpha > -\frac{1}{2}$, it follows that $I_{1,1}(y) \lesssim 1$. In conclusion,

$$
I_1(y) \lesssim 1.
\tag{3.11}
$$

Next we estimate $I_2(y)$. In $\Omega \setminus \Omega_{\epsilon_0}(F)$, $\varphi_{\gamma V} \sim \varphi_0$. Therefore by [14, Thm.1.5], for any $y \in \partial\Omega$ and $x \in \Omega \cap B_{r_0}(y) \setminus \Omega_{\epsilon_0}(F)$,

$$
\begin{aligned}
\frac{\varphi_{\gamma V}(x)}{\varphi_0(x)} K_{\gamma V}(x,y) &\sim \frac{\varphi_{\gamma V}(x)^2}{\varphi_0(x)\varphi_{\gamma V}(x_y)^2} |x-y|^{2-N} \\
&\sim \frac{\delta(x)}{|x-y|^2} |x-y|^{2-N} \sim K_0(x,y).
\end{aligned}
\tag{3.12}
$$

It follows that

$$
\begin{aligned}
I_2(y) &= \int_{\Sigma_\beta \cap B_{r_0}(y) \setminus \Omega_{\epsilon_0}(F)} \frac{\varphi_{\gamma V}(x)}{\varphi_0(x)} K_{\gamma V}(x,y) dS_x \\
&\quad + \int_{(\Sigma_\beta \setminus B_{r_0}(y)) \setminus \Omega_{\epsilon_0}(F)} \frac{\varphi_{\gamma V}(x)}{\varphi_0(x)} K_{\gamma V}(x,y) dS_x \\
&\sim \int_{\Sigma_\beta \cap B_{r_0}(y) \setminus \Omega_{\epsilon_0}(F)} \frac{\delta(x)}{|x-y|^2} |x-y|^{2-N} dS_x \\
&\quad + \int_{(\Sigma_\beta \setminus B_{r_0}(y)) \setminus \Omega_{\epsilon_0}(F)} K_0(x,y) dS_x \\
&\sim \int_{\Sigma_\beta \setminus \Omega_{\epsilon_0}(F)} K_0(x,y) dS_x \lesssim 1.
\end{aligned}
\tag{3.13}
$$

This fact and (3.11) imply that

$$
I(y) \lesssim 1.
\tag{3.14}
$$

On the other hand, from (3.8) and (3.13),

$$
I(y) = I_1(y) + I_2(y) \gtrsim \int_{\Sigma_\beta \cap B_{\epsilon_0}(y)} K_0(x,y) dS_x \gtrsim 1.
\tag{3.15}
$$

Thus $I(y) \sim 1$.

*Case 2:* $\alpha > 0$. By combining (3.6), (3.7) and (3.13), we obtain

$$
\begin{aligned}
I(y) &= I_1(y) + I_2(y) \\
&\lesssim \int_{\Sigma_{\beta,\epsilon_0}(F) \cap B_{\epsilon_0}(y)} K_0(x,y) dS_x + \epsilon_0^{1+2\alpha} + \int_{\Sigma_\beta \setminus \Omega_{\epsilon_0}(F)} K_0(x,y) dS_x \\
&\lesssim 1.
\end{aligned}
\tag{3.16}
$$

By (3.6) and the inequality $\delta(x) \leq \delta_F(x)$ we obtain

$$
\begin{aligned}
I_1(y) &\gtrsim \int_{\Sigma_{\beta,\epsilon_0}(F) \cap B_{r_0}(y)} \left( \frac{\delta(x)}{\delta(x)+|x-y|} \right)^{2\alpha} \delta(x) |x-y|^{-N} dS_x \\
&\gtrsim \int_{\Sigma_{\beta,\epsilon_0}(F) \cap B_{2\beta}(y)} \left( \frac{\delta(x)}{\delta(x)+|x-y|} \right)^{2\alpha} \delta(x) |x-y|^{-N} dS_x \\
&\gtrsim \int_{\Sigma_{\beta,\epsilon_0}(F) \cap B_{2\beta}(y)} K_0(x,y) dS_x.
\end{aligned}
\tag{3.17}
$$

Hence $I(y) = I_1(y) + I_2(y) \gtrsim 1$. Consequently, by (3.16),

$$
I(y) \sim 1.
$$

$\square$

**Proof of Theorem 1.1.** The theorem follows from Theorem 3.1 and Fubini theorem. $\square$

**Proposition 3.2.** *There exists a positive constant* $c = c(N, \gamma, F, \Omega)$ *such that*

$$(3.18) \qquad \int_{\Sigma_\beta} \frac{\varphi_{\gamma V}(x)}{\varphi_0(x)} \mathbb{G}_{\gamma V}[\tau](x) dS_x \leq c \int_\Omega \varphi_{\gamma V} d\tau \quad \forall \beta \in (0, r_0)$$

*and*

$$(3.19) \qquad \operatorname{tr}^*(\mathbb{G}_{\gamma V}[\tau]) = 0 \qquad \forall \tau \in \mathfrak{M}^+(\Omega, \varphi_{\gamma V}).$$

*Proof. Case 1:* $\alpha \leq 0$.

*Proof of* (3.18). By Fubini's theorem,

$$I := \int_{\Sigma_\beta} \frac{\varphi_{\gamma V}(x)}{\varphi_0(x)} \mathbb{G}_{\gamma V}[\tau](x) dS_x = \int_\Omega \int_{\Sigma_\beta} \frac{\varphi_{\gamma V}(x)}{\varphi_0(x)} G_{\gamma V}(x, y) dS_x d\tau(y).$$

Let $\hat{c} = 16(1 + \kappa)^2$ and denote:

$$E_\beta(y) := \Sigma_\beta \cap B_{\hat{c}\beta}(y),$$
$$E_{\beta,1}(y) := \{x \in E_\beta(y) : |x - y| \leq \hat{c}\delta(y)\},$$
$$E_{\beta,2}(y) := \{x \in E_\beta(y) : |x - y| \geq \hat{c}\delta(y)\} \quad \forall y \in \Omega.$$

Put

$$I_{\beta,1} := \int_\Omega \int_{E_{\beta,1}(y)} \frac{\varphi_{\gamma V}(x)}{\varphi_0(x)} G_{\gamma V}(x, y) dS_x d\tau(y),$$
$$I_{\beta,2} := \int_\Omega \int_{E_{\beta,2}(y)} \frac{\varphi_{\gamma V}(x)}{\varphi_0(x)} G_{\gamma V}(x, y) dS_x d\tau(y).$$

Let $y \in \Omega$ and $x \in E_{\beta,1}(y)$. Then

$$(3.20) \qquad \frac{1}{\hat{c} + 1}\beta \leq \delta(y) \leq (\hat{c} + 1)\beta$$

Indeed, for $x \in E_\beta(y)$,

$$\delta(y) \leq \delta(x) + |x - y| \leq \beta + \hat{c}\beta = (\hat{c} + 1)\beta.$$

Moreover, either $\beta \leq \delta(y)$ or

$$0 < \beta - \delta(y) \leq |x - y| \leq \hat{c}\delta(y) \Rightarrow \beta \leq (1 + \hat{c})\delta(y).$$

Similarly, if $x \in E_{\beta,2}(y)$ then

$$\hat{c}\beta > |x - y| \geq \hat{c}\delta(y) \Rightarrow \delta(y) < \beta$$

and consequently

$$(3.21) \qquad |x - y| \geq \beta - \delta(y) \geq \beta - \frac{1}{\hat{c}}|x - y| \Rightarrow \beta \leq (1 + \frac{1}{\hat{c}})|x - y|.$$

It follows that if $x \in E_{\beta,1}(y)$ then there exists a constant $c_1 > 0$, independent of $\beta$, such that,

$$(3.22) \qquad \frac{1}{c_1}\varphi_{\gamma V}(y) \leq \varphi_{\gamma V}(x) \leq c_1 \varphi_{\gamma V}(y) \quad \forall x \in E_{\beta,1}(y).$$

This is proved by a Harnack chain argument: $x \in E_{\beta,1}(y)$ implies (3.20). Thus the length of the segment $[x, y]$ is at most $\hat{c}\beta$ and the segment is at distance at least $\beta/(1+\hat{c})$ form $\partial\Omega$. Consequently $[x, y]$ can be covered by a Harnack chain of $m$ balls of radius $\beta/2(1+\hat{c})$, such that $m$ depends on $\hat{c}$ but is independent of $\beta$. This implies (3.22).

*Estimate of $I_{\beta,1}$.* By [14, Theorem 1.3],

$$\int_{E_{\beta,1}(y)} G_{\gamma V}(x, y) dS_x \sim \int_{E_{\beta,1}(y)} |x - y|^{2-N} dS_x$$
$$\sim \int_0^{\hat{c}\beta} t^{2-N} t^{N-2} dt = \hat{c}\beta \quad \forall y \in \Omega.$$

This relation and (3.22) imply,

$$(3.23) \qquad I_{\beta,1} \lesssim \int_{\Omega_{(1+\hat{c})\beta}} \varphi_{\gamma V}(y) \, d\tau(y).$$

*Estimate of $I_{\beta,2}$.* We split the integral into two parts:

$$I_{\beta,2,1} := \int_\Omega \int_{E_{\beta,2}(y) \cap \Omega_{2\epsilon_0}(F)} \frac{\varphi_{\gamma V}(x)}{\varphi_0(x)} G_{\gamma V}(x, y) dS_x d\tau(y),$$
$$I_{\beta,2,2} := \int_\Omega \int_{E_{\beta,2}(y) \setminus \Omega_{2\epsilon_0}(F)} \frac{\varphi_{\gamma V}(x)}{\varphi_0(x)} G_{\gamma V}(x, y) dS_x d\tau(y).$$

If $x \in E_{\beta,2}(y)$ then there exists $b_0 \in (0, 1)$ such that

$$\min(\delta(x), \delta(y)) \leq b_0 |x - y|$$

and therefore, by [14, Theorem 1.4],

$$G_{\gamma V}(x, y) \sim |x - y|^{2-N} \frac{\varphi_{\gamma V}(x) \varphi_{\gamma V}(y)}{\varphi_{\gamma V}(\tilde{x}) \varphi_{\gamma V}(\tilde{y})}$$

where $\tilde{x}, \tilde{y}$ are arbitrary points in $A(x, y)$. We choose $\tilde{x} = \tilde{y}$ such that

$$(3.24) \qquad \delta(\tilde{x}) \sim |x - y|, \quad \delta_F(\tilde{x}) \lesssim \delta_F(x) + |x - y|.$$

If $x \in E_{\beta,2}(y) \cap \Omega_{2\epsilon_0}(F)$ then by Lemma 2.2 and the fact that $\alpha < 0$,
$$(3.25)$$
$$\frac{\varphi_{\gamma V}(x)}{\varphi_0(x)} G_{\gamma V}(x, y) \sim \delta_F(x)^\alpha |x - y|^{2-N} \frac{\delta(x) \delta_F(x)^\alpha}{\delta(\tilde{x})^2 \delta_F(\tilde{x})^{2\alpha}} \varphi_{\gamma V}(y)$$
$$\sim \delta(x) |x - y|^{-N} \left( \frac{\delta_F(x)}{\delta_F(\tilde{x})} \right)^{2\alpha} \varphi_{\gamma V}(y)$$
$$\lesssim \delta(x) |x - y|^{-N} \left( \frac{\delta_F(x)}{\delta_F(x) + |x - y|} \right)^{2\alpha} \varphi_{\gamma V}(y).$$

Finally, since $\delta(x) \leq \delta_F(x)$ and (by assumption) $\alpha < 0$,

$$(3.26) \qquad \frac{\varphi_{\gamma V}(x)}{\varphi_0(x)} G_{\gamma V}(x, y) \lesssim \delta(x) |x - y|^{-N} \left( \frac{\delta(x)}{\delta(x) + |x - y|} \right)^{2\alpha} \varphi_{\gamma V}(y)$$
$$\lesssim \beta |x - y|^{-N} (1 + (|x - y|/\beta)^{-2\alpha}) \varphi_{\gamma V}(y).$$

Using (3.21) and the fact that $|x - y| < \hat{c}\beta$

$$(3.27) \qquad \int_{E_{\beta,2} \cap \Omega_{\epsilon_0}} \beta |x - y|^{-N} dS_x \lesssim \int_{\beta\hat{c}/(\hat{c}+1)}^{\hat{c}\beta} \beta r^{-N} r^{N-2} dr \lesssim 1.$$

This inequality and (3.26) imply,

$$(3.28) \qquad I_{\beta,2,1} \lesssim \int_{\Omega_{(1+\hat{c})\beta}} \varphi_{\gamma V}(y) \, d\tau(y).$$

In $\Omega \setminus \Omega_{2\epsilon_0}(F)$ we have $\varphi_{\gamma V} \sim \varphi_0$. If, in addition $|x - y| < \hat{c}\beta$, then $y \in \Omega \setminus \Omega_{\epsilon_0}(F)$. (We assume, as we may, that $\hat{c}\beta < \epsilon_0$.) Therefore if $x \in E_{\beta,2} \setminus \Omega_{2\epsilon_0}(F)$, there exists a constant $c_2 > 0$ such that

$$\frac{1}{c_2} G_{\gamma V}(x, y) \leq G_0(x, y) \leq c_2 G_{\gamma V}(x, y).$$

This and (3.27) imply

$$(3.29) \qquad \begin{aligned} I_{\beta,2,2} &\lesssim \int_\Omega \int_{E_{\beta,2}(y) \setminus \Omega_{2\epsilon_0}(F)} G_0(x, y) dS_x d\tau(y) \\ &\lesssim \int_\Omega \int_{E_{\beta,2}(y) \setminus \Omega_{2\epsilon_0}(F)} \beta\delta(y) |x - y|^{-N} dS_x d\tau(y) \\ &\lesssim \int_{\Omega_{(1+\hat{c})\beta}} \varphi_{\gamma V}(y) \, d\tau(y). \end{aligned}$$

By (3.23), (3.28) and (3.29) we conclude that,

$$(3.30) \qquad \begin{aligned} I_{\beta,1} + I_{\beta,2} &= \int_\Omega \int_{\Sigma_\beta \cap B_{\hat{c}\beta}(y)} \frac{\varphi_{\gamma V}(x)}{\varphi_0(x)} G_{\gamma V}(x, y) dS_x d\tau(y) \\ &\lesssim \int_{\Omega_{(1+\hat{c})\beta}} \varphi_{\gamma V} d\tau \to 0 \quad \text{as } \beta \to 0. \end{aligned}$$

Next we estimate the integral $I'_\beta$ defined below.

$$I'_\beta := \int_\Omega \int_{E'_\beta(y)} \frac{\varphi_{\gamma V}(x)}{\varphi_0(x)} G_{\gamma V}(x, y) dS_x d\tau(y), \quad E'_\beta(y) := \Sigma_\beta \setminus B_{\hat{c}\beta}(y).$$

If $x \in E'_\beta(y)$ then the conditions of [14, Theorem 1.4] hold and consequently $G_{\gamma V}$ is estimated by

$$(3.31) \qquad G_{\gamma V}(x, y) \sim |x - y|^{2-N} \frac{\varphi_{\gamma V}(x)\varphi_{\gamma V}(y)}{\varphi_{\gamma V}(\tilde{x})\varphi_{\gamma V}(\tilde{y})}.$$

As before we choose $\tilde{x} = \tilde{y}$ such that (3.24) holds.

We estimate $I'_\beta$, splitting the domain of integration into three parts. Let

$$(3.32) \qquad I'_{\beta,1} := \int_{\Omega_{\epsilon_0}(F)} \int_{E'_\beta(y) \setminus \Omega_{2\epsilon_0}(F)} \frac{\varphi_{\gamma V}(x)}{\varphi_0(x)} G_{\gamma V}(x, y) dS_x d\tau(y)$$

Then $|x - y| \geq \epsilon_0$ and consequently $\delta(\tilde{x}) \gtrsim \epsilon_0$. This implies that there exists a constant $c(\epsilon_0)$ such that $\varphi_{\gamma V}(\tilde{x}) \geq c(\epsilon_0)$ and

$$G_{\gamma V}(x, y) \lesssim \varphi_{\gamma V}(x) \varphi_{\gamma V}(y).$$

Further, in $\Omega \setminus \Omega_{2\epsilon_0}(F)$, $\varphi_{\gamma V} \sim \varphi_0$ so that

$$\frac{\varphi_{\gamma V}(x)}{\varphi_0(x)} G_{\gamma V}(x, y) \lesssim \delta(x) \varphi_{\gamma V}(y).$$

We conclude that,

$$I'_{\beta,1} \lesssim \beta \int_\Omega \varphi_{\gamma V}(y) d\tau(y). \tag{3.33}$$

We turn to the estimate of

$$I'_{\beta,2} := \int_{\Omega_{\epsilon_0}(F)} \int_{E'_\beta(y) \cap \Omega_{2\epsilon_0}(F)} \frac{\varphi_{\gamma V}(x)}{\varphi_0(x)} G_{\gamma V}(x, y) dS_x d\tau(y). \tag{3.34}$$

As in the estimate of $I_{\beta,2}$ we obtain (3.26). Since $|x - y| > \hat{c}\beta$ it follows that

$$\frac{\varphi_{\gamma V}(x)}{\varphi_0(x)} G_{\gamma V}(x, y) \lesssim \beta^{1+2\alpha} |x - y|^{-N-2\alpha} \varphi_{\gamma V}(y).$$

Hence,

$$\begin{aligned} I'_{\beta,2} &\lesssim \int_{\Omega_{\epsilon_0}(F)} \int_{\Omega_{2\epsilon_0}(F) \setminus B_{\hat{c}\beta}(y)} \beta^{1+2\alpha} |x - y|^{-N-2\alpha} dS_x \varphi_{\gamma V}(y) d\tau(y) \\ &\lesssim \int_{\Omega_{\epsilon_0}(F)} \beta^{1+2\alpha} \int_{\hat{c}\beta}^\infty r^{-2-2\alpha} dr \varphi_{\gamma V}(y) d\tau(y). \end{aligned}$$

In conclusion,

$$I'_{\beta,2} \lesssim \int_{\Omega_{\epsilon_0}(F)} \varphi_{\gamma V}(y) d\tau(y). \tag{3.35}$$

It remains to estimate,

$$I'_{\beta,3} := \int_{\Omega \setminus \Omega_{\epsilon_0}(F)} \int_{E'_\beta(y)} \frac{\varphi_{\gamma V}(x)}{\varphi_0(x)} G_{\gamma V}(x, y) dS_x d\tau(y). \tag{3.36}$$

In $\Omega \setminus \Omega_{\frac{\epsilon_0}{2}}(F)$, $\varphi_{\gamma V} \sim \varphi_0$ and $G_{\gamma V} \sim G_0$. Therefore,

$$\int_{\Omega \setminus \Omega_{\epsilon_0}(F)} \int_{E'_\beta(y) \setminus \Omega_{\frac{\epsilon_0}{2}}(F)} \frac{\varphi_{\gamma V}(x)}{\varphi_0(x)} G_{\gamma V}(x,y) dS_x d\tau(y)$$

$$\sim \int_{\Omega \setminus \Omega_{\epsilon_0}(F)} \int_{E'_\beta(y) \setminus \Omega_{\frac{\epsilon_0}{2}}(F)} |x-y|^{-N} \varphi_0(x) \varphi_0(y) dS_x d\tau$$

$$\lesssim \int_\Omega \beta \int_{\hat{c}\beta}^\infty r^{-2} dr \, \varphi_{\gamma V} d\tau$$

$$\lesssim \int_\Omega \varphi_{\gamma V} d\tau.$$

For $y \in \Omega \setminus \Omega_{\epsilon_0}(F)$ and $x \in E'_\beta(y) \cap \Omega_{\frac{\epsilon_0}{2}}(F)$, $|x-y| \geq \frac{\epsilon_0}{2}$. Hence the conditions of [14, Theorem 1.4] hold and $G_{\gamma V}$ satisfies (3.31). This, together with Lemma 2.2, yields

$$\frac{\varphi_{\gamma V}(x)}{\varphi_0(x)} G_{\gamma V}(x,y) \sim \frac{\varphi_{\gamma V}}{\varphi_0} \frac{\varphi_{\gamma V}(x) \varphi_{\gamma V}(y)}{\varphi_{\gamma V}(\bar{x})^2} |x-y|^{2-N}$$

$$\lesssim \frac{\varphi_{\gamma V}(x)^2}{\varphi_0(x)} \varphi_{\gamma V}(y)$$

$$\sim \delta(x) \delta_F(x)^{2\alpha} \varphi_{\gamma V}(y).$$

Therefore

$$\int_{\Omega \setminus \Omega_{\epsilon_0}(F)} \int_{E'_\beta(y) \cap \Omega_{\frac{\epsilon_0}{2}}(F)} \frac{\varphi_{\gamma V}(x)}{\varphi_0(x)} G_{\gamma V}(x,y) \lesssim \beta^{1+2\alpha} \int_\Omega \varphi_{\gamma V}(y) d\tau(y).$$

Thus

$$(3.37) \qquad\qquad I'_{\beta,3} \lesssim \int_\Omega \varphi_{\gamma V}(y) d\tau(y).$$

Combining (3.30), (3.33), (3.35) and (3.37) we obtain (3.18).

*Proof of* (3.19). Given $\epsilon > 0$ choose $\beta_1 > 0$ such that

$$\int_{\Omega_{\beta_1}} \varphi_{\gamma V} d\tau < \epsilon.$$

Let $\tau_1 = \tau \chi_{\bar{\Omega}_{\beta_1}}$ and $\tau_2 = \tau - \tau_1$ and denote $v_i = \mathbb{G}_{\gamma V}[\tau_i]$, $i = 1, 2$.

By (3.18), there exists a constant $c(\beta_1)$ such that

$$(3.38) \qquad I^*_{\beta,1} := \int_{\Sigma_\beta} \frac{\varphi_{\gamma V}}{\varphi_0} v_1 dS \leq c\epsilon \quad \forall \beta \in (0, \epsilon_0/\hat{c}).$$

Next we estimate the integral

$$I^*_{\beta,2} := \int_{\Sigma_\beta} \frac{\varphi_{\gamma V}}{\varphi_0} v_2 dS$$

when $\beta < \beta_1/\hat{c}$. By [14, Theorem 1.4],

$$(3.39) \quad I^*_{\beta,2} \lesssim \int_{\Omega \setminus \Omega_{\beta_1}} \int_{\Sigma_\beta} \frac{\varphi_{\gamma V}(x)}{\varphi_0(x)} \frac{\varphi_{\gamma V}(x)\varphi_{\gamma V}(y)}{\varphi_{\gamma V}(\tilde{x})^2} |x-y|^{2-N} dS_x d\tau_2(y).$$

Here $|x - y| \geq \beta_1/2$ and $\delta(\tilde{x}) \sim \beta_1/2$. Therefore $\varphi_{\gamma V}(\tilde{x})^2 \geq a > 0$ for all $x, y$ as above. Consequently, by Lemma 2.2,

$$(3.40)$$
$$\int_{\Sigma_\beta \cap \Omega_{\epsilon_0}(F)} \frac{\varphi_{\gamma V}(x)}{\varphi_0(x)} \frac{\varphi_{\gamma V}(x)\varphi_{\gamma V}(y)}{\varphi_{\gamma V}(\tilde{x})^2} |x-y|^{2-N} dS_x \lesssim$$
$$\int_{\Sigma_\beta \cap \Omega_{\epsilon_0}(F)} \delta(x)\delta_F(x)^{2\alpha}\varphi_{\gamma V}(y) dS_x.$$

Since $\alpha \leq 0$, it follows that

$$(3.41) \quad \int_{\Sigma_\beta \cap \Omega_{\epsilon_0}(F)} \frac{\varphi_{\gamma V}(x)}{\varphi_0(x)} \frac{\varphi_{\gamma V}(x)\varphi_{\gamma V}(y)}{\varphi_{\gamma V}(\tilde{x})^2} |x-y|^{2-N} dS_x \lesssim \beta^{1+2\alpha}\varphi_{\gamma V}(y).$$

Outside $\Omega_{\epsilon_0}(F)$, $\varphi_{\gamma V} \sim \varphi_0$. Therefore,

$$(3.42) \quad \int_{\Sigma_\beta \setminus \Omega_{\epsilon_0}(F)} \frac{\varphi_{\gamma V}(x)}{\varphi_0(x)} \frac{\varphi_{\gamma V}(x)\varphi_{\gamma V}(y)}{\varphi_{\gamma V}(\tilde{x})^2} |x-y|^{2-N} dS_x \lesssim \beta\varphi_{\gamma V}(y).$$

By (3.39), (3.41), (3.42),

$$(3.43) \quad I^*_{\beta,2} \lesssim c_\beta \int_\Omega \varphi_{\gamma V}(y) d\tau(y)$$

where $c_\beta \to 0$ as $\beta \to 0$. Thus (3.19) follows from (3.38) and (3.43).

*Case 2:* $\alpha > 0$. First we observe that estimates (3.23) and (3.29) still hold with $\alpha > 0$. Next we estimate $I_{\beta,2,1}$.

If $x \in E_{\beta,2}(y)$ and $\tilde{x} \in A(x,y)$ then - see (3.1) -

$$\frac{1}{2}|x-y| \leq \delta(\tilde{x}) \leq \delta_F(\tilde{x}) \quad \text{and} \quad |\tilde{x} - \frac{x+y}{2}| \leq 2|x-y|.$$

It follows that

$$|\delta_F(\tilde{x}) - \delta_F(x)| \leq |\tilde{x} - x| \leq |\tilde{x} - \frac{x+y}{2}| + |x - \frac{x+y}{2}| \leq \frac{5}{2}|x-y|.$$

Hence

$$\delta_F(\tilde{x}) \geq \delta_F(x) - \frac{5}{2}|x-y| \geq \delta_F(x) - 5\delta_F(\tilde{x}) \Longrightarrow \delta_F(x) \leq 6\delta_F(\tilde{x}).$$

This inequality, the first part of (3.25) and the assumption $\alpha > 0$ lead to

$$I_{\beta,2,1} \lesssim \int_\Omega \int_{E_{\beta,2}(y) \cap \Omega_{2\epsilon_0}(F)} \delta(x)|x-y|^{-N}\varphi_{\gamma V}(y) dS_x d\tau(y).$$

By (3.27), we obtain (3.28) and therefore - since (3.29) holds independently of $\alpha$ - we obtain (3.30).

Next we estimate $I'_\beta$. We see that (3.33) and (3.37) are still valid with $\alpha > 0$. It remains to estimate $I'_{\beta,2}$. Since $\alpha > 0$, by the first part of (3.25), for any $y \in \Omega_{\epsilon_0}(F)$ and $x \in E'_\beta(y) \cap \Omega_{2\epsilon_0}(F)$,

$$\frac{\varphi_{\gamma V}(x)}{\varphi_0(x)} G_{\gamma V}(x,y) \lesssim \delta(x)|x-y|^{-N}\varphi_{\gamma V}(y).$$

Since $|x-y| > \hat{c}\beta$, it follows that

$$
\begin{aligned}
(3.44) \qquad I'_{\beta,2} &\lesssim \int_{\Omega_{\epsilon_0}(F)} \int_{\Omega_{2\epsilon_0}(F)\setminus B_{\hat{c}\beta}(y)} \beta|x-y|^{-N} dS_x \varphi_{\gamma V}(y) d\tau(y) \\
&\lesssim \int_{\Omega_{\epsilon_0}(F)} \beta \int_{\hat{c}\beta}^\infty r^{-2} dr \varphi_{\gamma V}(y) d\tau(y) \\
&\lesssim \int_{\Omega_{\epsilon_0}(F)} \varphi_{\gamma V}(y) d\tau(y).
\end{aligned}
$$

By combining (3.30), (3.33), (3.37) and (3.44), we obtain (3.18).

Finally, inequality (3.19) is proved in the same way as in Case 1. $\square$

**Corollary 3.3.** *There exists a positive constant $c$ depending on $N, \gamma, V$ and $\Omega$ such that*

$$(3.45) \quad \int_\Omega \frac{\varphi_{\gamma V}(x)}{\varphi_0(x)} \mathbb{G}_{\gamma V}[\tau](x) dx \leq c \int_\Omega \varphi_{\gamma V} d\tau \qquad \forall \tau \in \mathfrak{M}^+(\Omega, \varphi_{\gamma V}).$$

*Proof.* Since $\Omega$ is a bounded $C^2$ domain, this is an immediate consequence of Proposition 3.2. $\square$

## 4. Linear b.v.p.

In the sequel we assume, without further mention, that $\gamma$ satisfies (3.2). We start this section with a result concerning $L_{\gamma V}$ superharmonic function.

**Proposition 4.1.** *Let $u$ be a positive $L_{\gamma V}$ superharmonic function. Then there exist $\nu \in \mathfrak{M}^+(\partial\Omega)$ and $\tau \in \mathfrak{M}^+(\Omega, \varphi_{\gamma V})$ such that (1.14) and (1.15) hold. If $\operatorname{tr}^*(u) = 0$, $u$ is an $L_{\gamma V}$ potential (i.e. it does not dominate any positive $L_{\gamma V}$ harmonic function).*

*Proof.* By the Riesz decomposition theorem (see [1]) $u$ has a unique representation of the form

$$u = u_h + u_p$$

where $0 \leq u_h$ is $L_{\gamma V}$ harmonic and $0 \leq u_p$ is an $L_{\gamma V}$ potential. By the Martin representation theorem there exists $\nu \in \mathfrak{M}^+(\partial\Omega)$ such that $u_h = \mathbb{K}_{\gamma V}[\nu]$. Furthermore, it is known that every positive $L_{\gamma V}$ potential is the Green potential of a positive measure $\tau$. Since $\mathbb{G}_{\gamma V}[\tau]$ is finite, $\tau \in \mathfrak{M}^+(\Omega, \varphi_{\gamma V})$.

By Proposition 3.2 $\operatorname{tr}^*(u) = \nu$, which implies the last assertion. $\square$

**Proposition 4.2.** *For any $\tau \in \mathfrak{M}^+(\Omega, \varphi_{\gamma V})$ and $\nu \in \mathfrak{M}^+(\partial \Omega)$, problem* (1.11) *has a unique positive solution given by* (1.14).

*Proof.* Proposition 3.2 implies that the function $u$ given by (1.14) is a solution of (1.11).

If $u$ is a non-negative solution of (1.11) then $u$ is $L_{\gamma V}$ superharmonic and by Proposition 4.1, there exist $\tau' \in \mathfrak{M}^+(\Omega, \varphi_{\gamma V})$ and $\nu' \in \mathfrak{M}^+(\partial \Omega)$ such that

$$u = \mathbb{G}_{\gamma V}[\tau'] + \mathbb{K}_{\gamma V}[\nu'].$$

It follows that $\operatorname{tr}^*(u) = \nu'$ and $-L_{\gamma V} u = \tau'$. However by assumption, $\operatorname{tr}^*(u) = \nu$ and $-L_{\gamma V} u = \tau$. Hence $\nu' = \nu$ and $\tau' = \tau$ so that the solution $u$ satisfies (1.14). $\qquad \square$

**Remark 4.3.** As the problem is linear, the result can be extended in a standard way to problems with signed measures.

**Proposition 4.4.** *Let $w$ be a non-negative $L_{\gamma V}$ subharmonic function. If $w$ is dominated by an $L_{\gamma V}$ superharmonic function then*

$$L_{\gamma V} w =: \lambda \in \mathfrak{M}^+(\Omega, \varphi_{\gamma V})$$

*and $w$ has a normalized boundary trace $\nu \in \mathfrak{M}^+(\partial \Omega)$. Thus,*

$$(4.1) \qquad\qquad w + \mathbb{G}_{\gamma V}[\lambda] = \mathbb{K}_{\gamma V}[\nu].$$

*Proof.* The first assumption implies that there exists a positive Radon measure $\lambda$ in $\Omega$ such that $-L_{\gamma V} w = -\lambda$.

If $\lambda \in \mathfrak{M}^+(\Omega, \varphi_{\gamma V})$ then $v := w + \mathbb{G}_{\gamma V}[\lambda]$ is a non-negative $L_{\gamma V}$ harmonic function and consequently, by the representation theorem, $v = \mathbb{K}_{\gamma V}[\nu]$ for some $\nu \in \mathfrak{M}^+(\partial \Omega)$. By Proposition 3.2, $\operatorname{tr}^*(w) = \nu$ and (4.1) holds.

Next we verify that $\lambda \in \mathfrak{M}^+(\Omega, \varphi_{\gamma V})$. Let $\{\Omega_n\}$ be an exhaustion of $\Omega$ consisting of $C^2$ domains.

Put $\lambda_n = \lambda \mathbf{1}_{\Omega_n}$ and $h_n = w \lfloor_{\partial \Omega_n}$. Let $v_n$ be the unique solution of the boundary value problem,

$$(4.2) \qquad -L_{\gamma V} v = -\lambda_n \text{ in } \Omega_n, \quad v = h_n \text{ on } \partial \Omega_n.$$

As in Proposition 4.2,

$$(4.3) \qquad\qquad v_n = \mathbb{G}^{\Omega_n}_{\gamma V}[-\lambda_n] + \mathbb{P}^{\Omega_n}_{\gamma V}[h_n],$$

where $\mathbb{G}^{\Omega_n}_{\gamma V}$ is the Green operator and $\mathbb{P}^{\Omega_n}_{\gamma V}$ is the Poisson operator of $L_{\gamma V}$ in $\Omega_n$. Since $w$ is a solution of (4.2), $v_n = w$ in $\Omega_n$.

By assumption, there exists an $L_{\gamma V}$ superharmonic function, say $W$, such that $w \leq W$ in $\Omega$. Hence (see Proposition 4.1)

$$\mathbb{P}^{\Omega_n}_{\gamma V}[h_n] \leq \mathbb{P}^{\Omega_n}_{\gamma V}[W \lfloor_{\partial \Omega_n}] \leq W \quad \text{in } \Omega_n.$$

Therefore, by (4.3), $w + \mathbb{G}^{\Omega_n}_{\gamma V}[\lambda_n] \leq W$. Letting $n \to \infty$ we conclude that $\mathbb{G}^{\Omega}_{\gamma V}[\lambda] \leq W$, which implies that $\lambda \in \mathfrak{M}^+(\Omega, \varphi_{\gamma V})$. $\qquad \square$

**Proposition 4.5.** *Let $w$ be a positive $L_{\gamma V}$ subharmonic function. Then the following are equivalent:*

(i) *$w$ is dominated by an $L_{\gamma V}$ superharmonic function.*

(ii) *$w$ has a normalized boundary trace.*

(iii) *$L_{\gamma V} w = \tau \in \mathfrak{M}^+(\Omega, \varphi_{\gamma V})$.*

(iv) *$w$ is dominated by an $L_{\gamma V}$ harmonic function.*

*If (iii) holds, there exists a measure $\nu \in \mathfrak{M}^+(\partial\Omega)$ such that*

$$(4.4) \qquad w + \mathbb{G}_{\gamma V}[\tau] = \mathbb{K}_{\gamma V}[\nu].$$

*Thus,*

$$(4.5) \qquad \mathrm{tr}^*(w) = \nu \text{ and } w \le \mathbb{K}_{\gamma V}[\nu].$$

*Proof.* (i) $\Rightarrow$ (ii), (iii) by Proposition 4.4.

(iii) $\Rightarrow$ (4.4). We note that $w + \mathbb{G}_{\gamma V}[\tau]$ is positive $L_{\gamma V}$ harmonic and consequently there exists a measure $\nu \in \mathfrak{M}^+(\partial\Omega)$ such that (4.4) holds.

(iii) $\Rightarrow$ (ii). By (4.4) and Proposition 3.2, $\mathrm{tr}^*(w) = \nu$.

(iii) $\Rightarrow$ (iv). By (4.4) and since $\tau \ge 0$, $w \le \mathbb{K}_{\gamma V}[\nu]$.

(iv) $\Rightarrow$ (i). Obviously $\qquad\qquad\qquad\qquad\qquad\qquad\square$

**Proof of Theorem 1.4.** The theorem follows from Propositions 3.2, 4.4, 4.1 and 4.2. $\qquad\qquad\qquad\qquad\qquad\qquad\qquad\qquad\square$

## 5. Nonlinear problem: Moderate solutions and an a-priori estimate

Let $f \in C(\mathbb{R})$ be an odd, monotone increasing function. We recall that a function $u$ is a solution of the nonlinear equation (1.1) if $u \in L^1_{loc}(\Omega)$, $f(u) \in L^1_{loc}(\Omega)$ and the equation holds in the distribution sense. Evidently, if $u$ is a non-negative solution of (1.1) then $u$ is $L_{\gamma V}$ subharmonic. A positive solution of (1.1) is $L_{\gamma V}$ *moderate* if it is dominated by an $L_{\gamma V}$ harmonic function.

As in the previous section, it will be assumed, unless otherwise stated, that $\gamma$ satisfies (3.2).

**Proof of Theorem 1.6.** The theorem is an immediate consequence of Proposition 4.5. $\qquad\qquad\qquad\qquad\qquad\qquad\qquad\qquad\square$

For the proof of Theorem 1.7 we need the following.

**Lemma 5.1.** *Let $\tau \in \mathfrak{M}^+(\Omega, \varphi_{\gamma V})$. Then*

$$(5.1) \qquad \int_\Omega \mathbb{G}_{\gamma V}[\tau] \varphi_{\gamma V} dx = \frac{1}{\lambda_{\gamma V}} \int_\Omega \varphi_{\gamma V} d\tau$$

*where $\lambda_{\gamma V}$ is the first eigenvalue of $-L_{\gamma V}$.*

*Proof.* For $\beta > 0$, put

$$I(\beta) := -\lambda_{\gamma V} \int_{D_\beta} \mathbb{G}_{\gamma V}[\tau] \varphi_{\gamma V} \, dx + \int_{D_\beta} \varphi_{\gamma V} \, d\tau,$$

where $D_\beta = \{x \in \Omega : \delta(x) > \beta\}$.

The assumption on $\tau$ and Proposition 3.2 imply that,

$$(5.2) \qquad \lim_{\beta \to 0} \int_{\Omega_\beta} \varphi_{\gamma V} \, d\tau = 0, \quad \lim_{\beta \to 0} \int_{\Omega_\beta} \mathbb{G}_{\gamma V}[\tau] \varphi_{\gamma V} \, dx = 0.$$

Therefore it is sufficient to prove that, for every $\tau$ with compact support in $\Omega$,

$$(5.3) \qquad \lim_{\beta \to 0} I(\beta) = 0.$$

Suppose that $\operatorname{supp} \tau \subset D_{\bar\beta}$ and let $\beta \in (0, \bar\beta/2)$. Applying Green's theorem in $D_\beta$ we obtain

$$-\lambda_{\gamma V} \int_{D_\beta} \mathbb{G}_{\gamma V}[\tau] \varphi_{\gamma V} \, dx = \int_{D_\beta} \mathbb{G}_{\gamma V}[\tau] L_{\gamma V} \varphi_{\gamma V} \, dx$$

$$= -\int_{D_\beta} \varphi_{\gamma V} d\tau - \int_{\Sigma_\beta} \frac{\partial \mathbb{G}_{\gamma V}[\tau]}{\partial \mathbf{n}} \varphi_{\gamma V} dS_x + \int_{\Sigma_\beta} \frac{\partial \varphi_{\gamma V}}{\partial \mathbf{n}} \mathbb{G}_{\gamma V}[\tau] dS_x.$$

Thus

$$(5.4) \qquad I(\beta) = -\int_{\Sigma_\beta} \frac{\partial \mathbb{G}_{\gamma V}[\tau]}{\partial \mathbf{n}} \varphi_{\gamma V} dS_x + \int_{\Sigma_\beta} \frac{\partial \varphi_{\gamma V}}{\partial \mathbf{n}} \mathbb{G}_{\gamma V}[\tau] dS_x.$$

By interior elliptic estimates, for every $x \in \Sigma_\beta$,

$$\left| \frac{\partial \varphi_{\gamma V}}{\partial \mathbf{n}}(x) \right| \leq C \sup_{|\xi - x| < \beta/4} \varphi_{\gamma V}(\xi) \beta^{-1}.$$

Therefore by Harnack's inequality,

$$\left| \frac{\partial \varphi_{\gamma V}}{\partial \mathbf{n}}(x) \right| \leq C \varphi_{\gamma V}(x) \beta^{-1} \quad \forall x \in \Sigma_\beta.$$

Hence, by Proposition 3.2,

$$(5.5) \qquad \lim_{\beta \to 0} \int_{\Sigma_\beta} \frac{\partial \varphi_{\gamma V}}{\partial \mathbf{n}} \mathbb{G}_{\gamma V}[\tau] dS_x = 0.$$

Since $\operatorname{supp} \tau \subset D_{\bar\beta}$, the same argument as above yields,

$$\left| \frac{\partial \mathbb{G}_{\gamma V}[\tau]}{\partial \mathbf{n}}(x) \right| \leq C \mathbb{G}_{\gamma V}[\tau](x) \beta^{-1} \quad \forall x \in \Sigma_\beta.$$

Therefore applying again Proposition 3.2 we obtain

$$(5.6) \qquad \lim_{\beta \to 0} \int_{\Sigma_\beta} \frac{\partial \mathbb{G}_{\gamma V}[\tau]}{\partial \mathbf{n}} \varphi_{\gamma V} dS_x = 0.$$

Finally, (5.4) – (5.6) imply (5.3).                    □

**Proof of Theorem 1.7.**

(i) (Uniqueness) Let $u_1$ and $u_2$ be two positive solutions of (1.16). Then $v := (u_1 - u_2)_+$ is a nonnegative $L_{\gamma V}$ subharmonic function and $\operatorname{tr}^*(v) = 0$. Furthermore, by Theorem 1.6, $f(u_1), f(u_2) \in L^1(\Omega, \varphi_{\gamma V})$ and $v \le \mathbb{G}_{\gamma V}[f(u_1) + f(u_2)] =: \bar{v}$. The inequality follows from the fact that, by Proposition 4.2,

$$u_i^* := \mathbb{K}_{\gamma V}[\nu] - u_i = \mathbb{G}_{\gamma V}[f(u_i)] \ge 0$$

and $v = (u_2^* - u_1^*)_+$.

Obviously $\bar{v}$ is $L_{\gamma V}$ superharmonic. Since $\operatorname{tr}^*(v) = 0$, Proposition 4.4 implies that $v = 0$. Thus $u_1 \le u_2$ and similarly $u_2 \le u_1$.

(ii) (Monotonicity) As before, $v := (u_1 - u_2)_+$ is $L_{\gamma V}$ subharmonic and it is dominated by an $L_{\gamma V}$ superharmonic function. Since $\nu_1 \le \nu_2$, $\operatorname{tr}^*(v) = 0$. Hence by Proposition 4.4, $v = 0$.

(iii) (A-priori estimate) Let $u$ be a positive solution of (1.16). By definition, $u$ has a normalized boundary trace. Therefore, by Theorem 1.6, $f(u) \in L^1(\Omega, \varphi_{\gamma V})$ and (1.17) holds.

Inequality (1.18) follows from (1.17) and (**??**).

The left hand side of (1.19) follows from (1.18), (**??**) and (3.45). The right hand side of (1.19) follows from (1.18) and Lemma 5.1. □

**Lemma 5.2.** *Assume $\gamma < C_H(V)$ and $D \Subset \Omega$ is a $C^2$ domain. Let $u_1, u_2 \in H^1_{loc}(D) \cap C(D)$ be respectively a positive subsolution and a positive supersolution of (1.1) in $D$ such that*

$$\limsup_{x \to \partial D}(u_1(x) - u_2(x)) = 0.$$

*Then $u_1 \le u_2$ a.e. in $D$.*

*Proof.* This is proved as in [5, Lemma 3.2]. □

**Lemma 5.3.** *Assume $\gamma < C_H(V)$ and $D \Subset \Omega$ is a $C^2$ domain. If $h \in C(\partial D)$ and $h \ge 0$, there exists a unique solution of the problem*

$$(5.7) \qquad \begin{cases} -L_{\gamma V} u + f(u) = 0 & in \ D \\ \qquad\qquad\quad u = h & on \ \partial D. \end{cases}$$

*Proof.* This is proved as in [16, Lemma 3.4]. □

**Proof of Proposition 1.9.** Let $u_0 := \mathbb{K}_{\gamma V}[\nu]$. Clearly $u_0 \in C^2(\Omega)$ and by assumption $f(u_0) \in L^1(\Omega, \varphi_{\gamma V})$.

Put $h_\beta := u_0 \lfloor_{\Sigma_\beta}$ where $\Sigma_\beta$ is given in (1.7) and. By Lemma 5.3 and Lemma 5.2 there exists a unique solution $u_\beta$ of (5.7) with $D$ replaced by $D_\beta$ and $h$ replaced by $h_\beta$, $\beta \in (0, \epsilon_0)$. Furthermore, by Lemma 5.2, as $u_0$ is a supersolution of (1.1), $\{u_\beta\}$ decreases as $\beta \downarrow 0$. Therefore $u := \lim_{\beta \to 0} u_\beta$ is a solution of (1.1) and

$$(5.8) \qquad u_\beta + \mathbb{G}_{\gamma V}^{D_\beta}[f(u_\beta)] = \mathbb{P}_{\gamma V}^{D_\beta}[h_\beta] = u_0$$

where $\mathbb{G}_{\gamma V}^{D_\beta}$ and $\mathbb{P}_{\gamma V}^{D_\beta}$ denote the Green and Poisson operators of $-L_{\gamma V}$ in $D_\beta$. Since $f(u_0) \in L^1(\Omega, \varphi_{\gamma V})$ and $u_\beta < u_0$ the sequence $\{\mathbb{G}_{\gamma V}^{D_\beta}[f(u_\beta)]\}$ is bounded. Letting $\beta \to 0$ in (5.8), we obtain

$$(5.9) \qquad u + \mathbb{G}_{\gamma V}[f(u)] = u_0 := \mathbb{K}_{\gamma V}[\nu].$$

$\square$

## 6. Subcritical nonlinearities

Let $\mathcal{F}(a^*)$ be defined as in Definition 1.12.

**Theorem 6.1.** *If $f \in \mathcal{F}(a^*)$, there exist positive constants $a_0, a_1$ such that, every positive solution $u$ of* (1.1) *in $\Omega$ satisfies the inequality,*

$$(6.1) \qquad f(u(x)/a_1) \leq a_0 \delta(x)^{-2} u(x) \quad \forall x \in \Omega.$$

*Proof.* Put $h(t) := f(t)/t$. By assumption $h$ is non-decreasing on $(0, \infty)$. If $u$ is a positive solution of (1.1) then

$$-\Delta u - \frac{|\gamma|}{\delta^2} u + h(u)u \leq 0.$$

By Harnack's inequality there exists $c_1 > 0$ such that, for every $y \in \Omega$,

$$(6.2) \qquad \sup_{D^y} u(x) \leq c_1 \inf_{D^y} u(x), \quad x \in D^y := B_{\delta(y)/2}(y).$$

Therefore, for every $y \in \Omega$,

$$(6.3) \qquad -\Delta u(x) - \frac{4|\gamma|}{\delta(y)^2} u(x) + h(u(y)/c_1)u(x) \leq 0 \quad \forall x \in D^y.$$

If

$$(6.4) \qquad \frac{4a|\gamma|}{\delta(y)^2} < \frac{1}{2} h(u(y)/c_1)$$

then

$$-\Delta u(x) + \frac{1}{2} h(u(y)/c_1)u(x) \leq 0 \quad \forall x \in D^y.$$

Therefore, applying again (6.2),

$$-\Delta u(x) + \frac{1}{2} h(u(x)/c_1^2)u(x) \leq 0 \quad \forall x \in D^y.$$

Hence, if $v(\xi) = u(\sqrt{2}\xi + y)/c_1^2$ we obtain,

$$(6.5) \qquad -\Delta v + f(v) \leq 0, \quad |\xi| < \frac{1}{2\sqrt{2}} \delta(y) \quad \forall y \in \Omega.$$

Thus $v$ satisfies (1.23) in the ball $|\xi| < \frac{1}{2\sqrt{2}} \delta(y)$. In particular, for $\xi = 0$,

$$f(v(0)) \leq 8a^* \delta(y)^{-2} v(0),$$

which translates to

$$(6.6) \qquad f(u(y)/c_1^2) \leq \frac{8a^*}{c_1^2} \delta(y)^{-2} u(y).$$

This inequality holds at every point $y \in \Omega$ such that (6.4) holds. At any other point $y \in \Omega$ we have,

$$(6.7) \qquad f(u(y)/c_1) \leq \frac{8|\gamma|}{c_1} \delta(y)^{-2} u(y).$$

<div style="text-align: right">□</div>

We recall that $\mathcal{M}_{\gamma V}(f)$ denotes the family of good measures (see Definition 1.8). The following is an immediate consequence of Theorem 6.1:

**Corollary 6.2.** *Assume $f \in \mathcal{F}(a^*)$. Let $\{\nu_n\} \subset \mathcal{M}_{\gamma V}(f)$ be an increasing sequence of measures. Denote by $u_n$ the solution of* (1.16) *with $\nu$ replaced by $\nu_n$. Then $u := \lim u_n$ is a solution of* (1.1)*. If the sequence $\{\nu_n\}$ is uniformly bounded then $\nu = \lim \nu_n \in \mathcal{M}_{\gamma V}(f)$ and* $\mathrm{tr}^*(u) = \nu$.

*Proof.* As $\{\nu_n\}$ is increasing, the sequence $\{u_n\}$ is increasing. By (6.1), $\{u_n\}$ is bounded in every compact subset of $\Omega$. (Recall that $f(t)/t$ is non-decreasing and tends to infinity as $t \to \infty$). Therefore the sequence converges and $u = \lim u_n$ is a solution of (1.1).

If, in addition, $\{\nu_n\}$ is bounded and $\nu = \lim \nu_n$ then $\mathbb{K}_{\gamma V}[\nu_n] \to \mathbb{K}_{\gamma V}[\nu]$. By Theorem 1.6,

$$u_n + \mathbb{G}_{\gamma V}[f(u_n)] = \mathbb{K}_{\gamma V}[\nu_n].$$

Hence,

$$u + \mathbb{G}_{\gamma V}[f(u)] = \mathbb{K}_{\gamma V}[\nu].$$

Thus $u$ is the solution of (1.16) and $\nu \in \mathcal{M}_{\gamma V}(f)$.

<div style="text-align: right">□</div>

**Proof of Theorem 1.13.** Put

$$E_n := \{y \in SC_{\gamma V}(f) : \|f(\mathbb{K}_{\gamma V}(\cdot, y))\|_{L^1(\Omega, \varphi_{\gamma V})} \leq n\}$$

and $\nu_n := \nu \mathbf{1}_{E_n}$. Then $f(\mathbb{K}_{\gamma V}[\nu_n]) \in L^1(\Omega, \varphi_{\gamma V})$ and consequently, by Proposition 1.9, $\nu_n \in \mathcal{M}_{\gamma V}(f)$. Let $u_n$ denote the solution of (1.16) when $\nu = \nu_n$. By Corollary 6.2, $u = \lim u_n$ is a solution of (1.1) and $\mathrm{tr}^*(u) = \lim \nu_n$. The assumption that $\nu$ is concentrated on $SC_{\gamma V}(f)$ implies that $\lim \nu_n = \nu$.

<div style="text-align: right">□</div>

**Proof of Theorem 1.14.** Put $v := \mathbb{K}_{\gamma V}[\nu] - u$. Then $u$ is positive and $L_{\gamma V}$ superharmonic. In addition $\mathrm{tr}^*(v) = 0$. Consequently,

$$\lim_{x \to y \in \partial\Omega} v(x)/\mathbb{K}_{\gamma V}[\nu](x) \to 0 \quad \text{non-tangentially } \nu\text{-a.e.}.$$

This implies (1.24).

<div style="text-align: right">□</div>

**Proof of Theorem 1.16.** Since $f(K_{\gamma V}(\cdot, y)) \in L^1(\Omega, \varphi_{\gamma V})$ and $f$ satisfies $\Delta_2$ condition, it follows that $f(kK_{\gamma V}(\cdot, y)) \in L^1(\Omega, \varphi_{\gamma V})$ for every $k > 0$. By Proposition 1.9, problem (1.21) has a (unique) solution, therefore $y \in SC_{\gamma V}(f)$.

Now assume that $y \in SC_{\gamma V}(f)$ and let $u$ be the solution of (1.21). By Theorem 1.14

$$(6.8) \qquad \lim_{x \to y} u(x)/K_{\gamma V}(x, y) = 1 \quad \text{non tangentially at } y.$$

Let $\mathcal{C}_y(a) = \{x \in \Omega : |x - y| < a\delta(x)\}$ where $a > 1$. Then there exists $r = r(a)$ such that

$$\mathcal{C}_y^r(a) := \mathcal{C}_y(a) \cap B_r(y) \subset \Omega \setminus \{y\} \Longrightarrow \overline{\mathcal{C}_y^r(a)} \setminus \{y\} \subset \Omega.$$

By Harnack's inequality and (6.8) there exists $C(a) > 0$ such that

$$(6.9) \qquad \mathbb{K}_{\gamma V}(x, y) \leq C(a)u(x) \quad \forall x \in \mathcal{C}_y^r(a).$$

Now suppose that $f(K_{\gamma V}(\cdot, y)) \notin L^1(\Omega, \varphi_{\gamma V})$. Then

$$\int_{\mathcal{C}_y^r(a)} \varphi_{\gamma V} f(C(a)u) \, dx \geq \int_{\mathcal{C}_y^r(a)} \varphi_{\gamma V} f(\mathbb{K}_{\gamma V}(x, y)) \, dx = \infty.$$

However, by Theorem 1.6, $f(u) \in L^1(\Omega, \varphi_{\gamma V})$. Since $f$ satisfies the $\Delta_2$ condition we reach a contradiction. $\qquad \square$

## 7. Existence and stability results

We recall that

$$q_F^* := \frac{N + 1 + \alpha}{N - 1 + \alpha} \quad \text{and} \quad q^* := \frac{N + 1}{N - 1}.$$

**Proposition 7.1.** *Assume $\gamma < C_H(V)$ and $q > 1$.*

*(i) If $y \in F$ and $-1 < \alpha$ then*

$$(7.1) \qquad K_{\gamma V}(\cdot, y) \in L^q(\Omega, \varphi_{\gamma V}) \Longleftrightarrow q < q_F^*.$$

*For every $q \in (1, q_F^*)$ there exists a constant $c(q)$ such that*

$$(7.2) \qquad \frac{1}{c(q)} \leq \|K_{\gamma V}(\cdot, y)\|_{L^q(\Omega, \varphi_{\gamma V})} \leq c(q) \quad \forall y \in F.$$

*(ii) If $y \in \partial\Omega \setminus F$ then*

$$(7.3) \qquad K_{\gamma V}(\cdot, y) \in L^q(\Omega, \varphi_{\gamma V}) \Longleftrightarrow q < q^*.$$

*Proof. Part (i):* $y \in F$. Let $x \in \Omega \cap B_{\epsilon_0}(y)$ and choose $x_y \in A(x, y)$ such that (3.4) holds. By [14, Theorem 1.5], Lemma 2.2 and (3.4),

$$(7.4) \qquad \begin{aligned} \varphi_{\gamma V}(x)K_{\gamma V}(x, y)^q &\sim \frac{\varphi_{\gamma V}(x)^{1+q}}{\varphi_{\gamma V}(x_y)^{2q}}|x - y|^{(2-N)q} \\ &\sim \frac{\delta(x)^{1+q}\delta_F(x)^{(1+q)\alpha}}{\delta(x_y)^{2q}\delta_F(x_y)^{2q\alpha}}|x - y|^{(2-N)q} \\ &\sim \frac{\delta(x)^{1+q}\delta_F(x)^{\alpha(1+q)}}{(\delta_F(x) + |x - y|)^{2q\alpha}}|x - y|^{-Nq} \end{aligned}$$

Since $y \in F$, $\delta_F(x) \leq |x-y|$, hence $\delta_F(x) + |x-y| \sim |x-y|$. Therefore

$$(7.5) \qquad \varphi_{\gamma V}(x) K_{\gamma V}(x,y)^q \sim \delta(x)^{1+q} \delta_F(x)^{\alpha(1+q)} |x-y|^{-(N+2\alpha)q}.$$

As $\delta \leq \delta_F \leq |x-y|$ it follows that:

- If $\alpha \leq 0$,

$$(7.6) \qquad \begin{aligned} \varphi_{\gamma V}(x) K_{\gamma V}(x,y)^q &\gtrsim \delta(x)^{1+q} |x-y|^{\alpha(1+q)-(N+2\alpha)q} \\ &= \left( \frac{\delta(x)}{|x-y|} \right)^{1+q} |x-y|^{1+\alpha-(N-1+\alpha)q}. \end{aligned}$$

- If $\alpha > 0$,

$$(7.7) \qquad \begin{aligned} \varphi_{\gamma V}(x) K_{\gamma V}(x,y)^q &\gtrsim \delta(x)^{(1+\alpha)(1+q)} |x-y|^{-(N+2\alpha)q} \\ &= \left( \frac{\delta(x)}{|x-y|} \right)^{(1+\alpha)(1+q)} |x-y|^{1+\alpha-(N-1+\alpha)q}. \end{aligned}$$

We also obtain the following estimates from above:

- If $-1 < \alpha \leq 0$,

$$(7.8) \qquad \begin{aligned} \varphi_{\gamma V}(x) K_{\gamma V}(x,y)^q &\lesssim \delta(x)^{(1+\alpha)(1+q)} |x-y|^{-(N+2\alpha)q} \\ &\leq |x-y|^{1+\alpha-(N-1+\alpha)q}. \end{aligned}$$

- If $\alpha > 0$,

$$(7.9) \qquad \begin{aligned} \varphi_{\gamma V}(x) K_{\gamma V}(x,y)^q &\lesssim \delta(x)^{1+q} |x-y|^{\alpha-(N+\alpha)q} \\ &\leq |x-y|^{1+\alpha-(N-1+\alpha)q}. \end{aligned}$$

Combining inequalities $(7.6) - (7.9)$, we conclude that for every $y \in F$ and $\alpha > -1$

$$(7.10) \qquad K_{\gamma V}(\cdot, y) \in L^q(\Omega, \varphi_{\gamma V}) \iff \int_0^1 t^{N+\alpha-q(N-1+\alpha)} dt < \infty.$$

This implies $(7.1)$ while inequalities $(7.6) - (7.9)$ imply $(7.2)$.

*Part (ii): $y \in \partial\Omega \setminus F$.* In this case there is a neighborhood $Q$ of $y$ such that in $Q \cap \Omega$, $\varphi_{\gamma V} \sim \varphi_0$ and $K_{\gamma V}(\cdot, y) \sim K_0(\cdot, y)$. Therefore the critical exponent is the same as in the case $V = 0$, i.e., $q^*$.         □

**Corollary 7.2.** *(a) Assume $\alpha > -1$ and $q \in (1, q_F^*)$. There exists $C_q > 0$ such that, if $\nu \in \mathfrak{M}^+(\partial\Omega)$ and $\operatorname{supp}\nu \subset F$ then,*

$$(7.11) \qquad \frac{1}{C_q} \|\nu\|_{\mathfrak{M}(\partial\Omega)} \leq \int_\Omega \mathbb{K}_{\gamma V}^q[\nu] \varphi_{\gamma V} \, dx \leq C_q \|\nu\|_{\mathfrak{M}(\partial\Omega)}.$$

*(b) Assume $\alpha > -1$ and $q \in (1, q^*)$. Let $r > 0$ and put*

$$(\partial\Omega)^r := \{y \in \partial\Omega : \operatorname{dist}(y, F) \geq r\}.$$

*Then there exists a constant $C_q(r)$ such that, for every $\nu \in \mathfrak{M}^+(\partial\Omega)$ with $\operatorname{supp}\nu \subset (\partial\Omega)^r$, $(7.11)$ remains valid if $C_q$ is replaced by $C_q(r)$.*

*Proof.* Statement (a) follows from (7.2). Statement (b) follows from the following facts: For $r > 0$ sufficiently small, there exists a constant $c(r)$ such that for every $y \in (\partial\Omega)^r$,

$$(7.12) \qquad \frac{1}{c(r)} K_0(x,y) \leq K_{\gamma V}(x,y) \leq c(r) K_0(x,y) \quad \forall x \in B_{r/2}(y),$$

$$\frac{1}{c(r)} \delta(x) \leq \varphi_{\gamma V}(x) \leq c(r)\delta(x) \quad \text{if } \delta_F(x) > r.$$

$\square$

**Proof of Theorem 1.17.** Statements (a) and (b) are an immediate consequence of Theorem 1.16 and Corollary 7.2. $\square$

Denote by $L_w^p(\Omega, \phi)$, $1 < p < \infty$, $0 < \phi \in C(\Omega)$ the weak $L^p$ space with weight $\phi$. (See e.g. [19, subsection 2.3.1] for more details.) For the next theorem we need,

**Lemma 7.3.** *Let $\alpha > -1$. For every $y \in \partial\Omega$ put*

$$\Lambda_y(x) = |x - y|^{-N+1-\alpha}, \ \eta_y(x) = \delta(x)|x - y|^\alpha \quad \forall x \in \Omega.$$

*Then*

$$(7.13) \qquad \|\Lambda_y\|_{L_w^{\frac{N+1+\alpha}{N-1+\alpha}}(\Omega, \eta_y)} \leq c$$

*where $c$ is a constant depending on $N, \alpha$ and $\mathrm{diam}\,(\Omega)$ but independent of $y \in \partial\Omega$.*

This lemma is similar to Lemma 2.3.2 of [19], except that the weight function is somewhat different. In particular, here the weight function depends on $y$ as well. However, throughout the proof $y$ is a fixed parameter and one must only notice that the constant $c$ is independent of $y$. In view of this observation and the fact that, for every $y \in \partial\Omega$,

$$\eta_y(x) \leq |x - y|^{1+\alpha} \quad \forall x \in \Omega,$$

the argument in the proof of [19, Lemma 2.3.2] applies also to the present case. Therefore the proof is ommitted.

**Theorem 7.4.** *Assume $\alpha > -1$.*

*(i) If $y \in F$ then $K_{\gamma V}(\cdot, y) \in L_w^{q_F^*}(\Omega, \varphi_{\gamma V})$ and*

$$(7.14) \qquad \|K_{\gamma V}(\cdot, y)\|_{L_w^{q_F^*}(\Omega, \varphi_{\gamma V})} \leq c,$$

*where $c$ is a constant independent of $y$.*

*(ii) If $y \in \partial\Omega \setminus F$ then $K_{\gamma V}(\cdot, y) \in L_w^{q^*}(\Omega, \varphi_{\gamma V})$ and, for every $r > 0$ there exists a constant $C_r$ such that,*

$$(7.15) \qquad \|K_{\gamma V}(\cdot, y)\|_{L_w^{q^*}(\Omega, \varphi_{\gamma V})} \leq C_r, \quad \forall y \in (\partial\Omega)^r$$

*where $(\partial\Omega)^r$ is defined as in Corollary 7.2.*

*Proof.* (i) $y \in F$. By [14, Theorem 1.5],

$$K_{\gamma V}(x, y) \sim \frac{\varphi_{\gamma V}(x)}{\varphi_{\gamma V}(x_y)^2} |x - y|^{2-N} \quad \forall x \in B_{\epsilon_0}(y) \cap \Omega_{\frac{r_0}{10\kappa}},$$

where $x_y$ can be chosen so that $\delta(x_y) \sim |x - y|$ and (since $y \in F$) $\delta_F(x_y) \leq C|x - y|$ (see (3.4)). Hence

$$
\begin{aligned}
(7.16) \qquad K_{\gamma V}(x, y) &\leq C \frac{|x - y|^{1+\alpha}}{|x - y|^{2(1+\alpha)}} |x - y|^{2-N} \\
&= C|x - y|^{1-N-\alpha} \qquad \forall x \in B_{\epsilon_0}(y) \cap \Omega_{\frac{r_0}{10\kappa}}.
\end{aligned}
$$

In $\Omega \setminus B_{\epsilon_0}(y)$, $K_\gamma(\cdot, y)$ is bounded by a constant. Therefore

$$(7.17) \qquad K_{\gamma V}(x, y) \lesssim |x - y|^{1-N-\alpha} \quad \forall x \in \Omega.$$

Since $y \in F$, $\max(\delta(x), \delta_F(x)) \leq |x - y|$. Therefore,

$$(7.18) \qquad \varphi_{\gamma V}(x) \leq \begin{cases} \delta(x)^{1+\alpha} & \text{if } \alpha \leq 0 \\ \delta(x)|x - y|^\alpha & \text{if } \alpha > 0. \end{cases}$$

By (7.17) and (7.18), using [19, Lemma 2.3.2] – specifically inequality (2.3.12) – when $\alpha < 0$ and Lemma 7.3 when $\alpha > 0$ we obtain (7.14).

(ii) $y \in \partial\Omega \setminus F$. Here the result follows from (7.12). $\qquad\square$

**Corollary 7.5.** *(i) If $\nu \in \mathfrak{M}^+(\partial\Omega)$ and $\operatorname{supp} \nu \subset F$ then $\mathbb{K}_{\gamma V}[\nu] \in L_w^{q_F^*}(\Omega, \varphi_{\gamma V})$ and*

$$(7.19) \qquad \|\mathbb{K}_{\gamma V}[\nu]\|_{L_w^{q_F^*}(\Omega, \varphi_{\gamma V})} \leq c \, \|\nu\|_{\mathfrak{M}(\partial\Omega)}$$

*with $c$ as in (7.14).*

*(ii) If $\nu \in \mathfrak{M}^+(\partial\Omega)$ and $\operatorname{supp} \nu \subset (\partial\Omega)^r$ then $\mathbb{K}_{\gamma V}[\nu] \in L_w^{q^*}(\Omega, \varphi_{\gamma V})$ and*

$$(7.20) \qquad \|\mathbb{K}_{\gamma V}[\nu]\|_{L_w^{q^*}(\Omega, \varphi_{\gamma V})} \leq C_r \, \|\nu\|_{\mathfrak{M}(\partial\Omega)}$$

*with $C_r$ as in (7.15).*

*Proof.* By Fubini's theorem this is an immediate consequence of Theorem 7.4. $\qquad\square$

**Proof of Theorem 1.18.** (i) By Theorem 7.4 and (1.27), if $y \in F$ and $p = q_F^*$ then

$$(7.21) \qquad f(kK_{\gamma V}(\cdot, y)) \in L^1(\Omega, \varphi_{\gamma V}) \quad \forall k > 0.$$

By Proposition 1.9, $F \subset SC_{\gamma V}(f)$. From (1.27) and (7.19), we deduce that for every finite positive measure $\nu$ supported in $F$, $f(\mathbb{K}_{\gamma V}[\nu]) \in L^1(\Omega, \varphi_{\gamma V})$. Again, by Proposition 1.9 we obtain that $\nu$ is a good measure.

The last statement follows from (7.19) and [19, Theorem 2.3.4].

(ii) Note that the assumption supp $\nu \subset \partial\Omega \backslash F$ implies that $\nu$ vanishes in a neighborhood of $F$. Therefore, using (7.20), the proof is the same as above.                                                    $\square$

**Proof of Theorem 1.19.** (i) and (ii). In view of Theorem 6.1, there exists a subsequence $\{u_{n'}\}$ that converges locally uniformly in $\Omega$ to a solution $u'$ of (1.1).

The sequence $\{\nu_n\}$ is bounded. Therefore, by Theorem 1.18, the sequence $\{f(u_n)\}$, which is dominated by $\{f(\mathbb{K}_{\gamma V}[\nu_n])\}$, is uniformly absolutely continuous in $L^1(\Omega, \varphi_{\gamma V})$. This fact and the pointwise convergence imply that

$$f(u_{n'}) \to f(u') \quad \text{in} \ \ L^1(\Omega, \varphi_{\gamma V}).$$

Consequently, by Proposition 3.2,

$$\mathbb{G}_{\gamma V}[f(u_{n'})] \to \mathbb{G}_{\gamma V}[f(u')] \quad \text{in} \ \ L^1(\Omega, \varphi_{\gamma V}).$$

Similarly, since $u_{n'} \leq \mathbb{K}_{\gamma V}[\nu_{n'}]$ and by Corollary 7.5, we deduce that the sequence $\{u_{n'}\}$ is uniformly absolutely continuous in $L^1(\Omega, \varphi_{\gamma V})$. This, combined with the pointwise convergence, yields

$$u_{n'} \to u' \quad \text{in} \ \ L^1(\Omega, \varphi_{\gamma V}).$$

Clearly

$$\mathbb{K}_{\gamma V}[\nu_n] \to \mathbb{K}_{\gamma V}[\nu] \quad \text{pointwise in} \ \Omega.$$

By Theorem 1.6,

$$u_{n'} + \mathbb{G}_{\gamma V}[f(u_{n'})] = \mathbb{K}_{\gamma V}[\nu_{n'}] \quad \text{a.e. in} \ \Omega.$$

This and the above convergence results imply that,

$$u' + \mathbb{G}_{\gamma V}[f(u')] = \mathbb{K}_{\gamma V}[\nu] \quad \text{a.e. in} \ \Omega.$$

Thus tr$^*(u') = \nu$. This means that $u'$ is a solution of (1.16). By uniqueness (see Theorem 1.7), $u' = u$ and hence the whole sequence $\{u_n\}$ satisfies (1.29).                                                    $\square$

**Acknowledgements.** Part of this work was carried out during a visit of the second author at the Technion. P. T. Nguyen is grateful to the Technion for the support and hospitality. The authors wish to thank Yehuda Pinchover for helpful discussions. The second author also wishes to thank F. Mahmoudi for some interesting discussions.

DEPARTMENT OF MATHEMATICS, TECHNION, HAIFA 32000, ISRAEL
    *E-mail address*: `marcusm@math.technion.ac.il`

DEPARTMENT OF MATHEMATICS, MASARYK UNIVERSITY, BRNO, CZECH REPUBLIC
    *E-mail address*: `ptnguyen@math.muni.cz;   nguyenphuoctai.hcmup@gmail.com`